\documentclass[12pt,reqno]{amsart}
\usepackage{amsmath}%
\usepackage{amsfonts}%
\usepackage{amssymb}%
\usepackage[margin=.5in]{geometry}
\usepackage{subfig}
\usepackage{array,multirow}
\usepackage{bm}
\usepackage{ifpdf}
\ifpdf
\usepackage{graphicx}
\else
\usepackage[draft]{graphicx}
\fi
\usepackage{xcolor}
\usepackage{enumerate}
\usepackage[mathscr]{euscript}
\usepackage[T1]{fontenc}
\usepackage{caption}
\usepackage{booktabs,multirow}
\usepackage{enumerate}
\usepackage{listings}
\usepackage{caption}
\usepackage{caption}
\usepackage{threeparttable}
\usepackage{bm}
\usepackage{csquotes}
\usepackage{float}
\usepackage{hyperref}
\MakeOuterQuote{"}
\hypersetup{
    colorlinks=true,
    linkcolor=blue,
    filecolor=magenta,      
    urlcolor=cyan,
		citecolor=red,
    pdfpagemode=UseNone,pdfstartview=
    }		
\definecolor{dkgreen}{rgb}{0,0.6,0}
\definecolor{gray}{rgb}{0.5,0.5,0.5}
\definecolor{mauve}{rgb}{0.58,0,0.82}

\AtBeginDocument{}
\AtBeginDocument{%
  \setlength\abovedisplayskip{20pt}
  \setlength\belowdisplayskip{20pt}}
\lstdefinestyle{myMathematica}{
frame=tb,
  language=Mathematica,
  aboveskip=3mm,
  belowskip=3mm,
  showstringspaces=false,
  columns=flexible,
  basicstyle={\ttfamily},
  numbers=none,
  numberstyle=\tiny\color{gray},
  keywordstyle=\color{blue},
  commentstyle=\color{mauve},
  stringstyle=\color{dkgreen},
  breaklines=true,
  breakatwhitespace=true,
	morekeywords={AsymptoticSolve},
  tabsize=3
 }
\newcommand\mscriptsize[1]{\mbox{\scriptsize\ensuremath{#1}}}

\begin{document}
\title[]{A numeric study of power expansions around singular points of algebraic functions, their radii of convergence, and accuracy profiles}

\author{Dominic C. Milioto}
\email{icorone@hotmail.com}
\date{\today}
\subjclass[2010]{Primary 1401; Secondary 1404} %
\begin{abstract}
An efficient method of computing power expansions of algebraic functions is the method of Kung and Traub \cite{Kung} and is based on exact arithmetic.  This paper shows a numeric approach is both feasible and accurate while also introducing a performance improvement to Kung and Traub's method based on the ramification extent of the expansions.  A new method is then described for computing radii of convergence using a series comparison test. Series accuracies are then fitted to a simple log-linear function in their domain of convergence and found to have low variance.   Algebraic functions up to degree $50$ were analyzed and timed.  A consequence of this work provided a simple method of computing the Riemann surface genus and was used as a cycle check-sum.  Mathematica ver. $13.2$ was used to acquire and analyze the data on a $4.0$ GHz quad-core desktop computer.
\end{abstract}
\maketitle
\section{Introduction}
The objective of this paper is five-fold:
\begin{enumerate}
  \item Provide a simple method of analyzing the branching geometry of algebraic functions,
		\item Describe a new method of determining radii of convergence,
	\item Construct accuracy and order functions of series expansions around singular points,
	\item Analyze test cases and summarize convergence, accuracy, and timing results,
	\item Provide an accessible teaching aid of this subject to interested readers.
	\end{enumerate}

The functions studied in this paper are algebraic functions $w(z)$ defined implicitly by the irreducible $n$-degree expression in $w$:
\begin{equation}
f(z,w)=a_0(z)+a_1(z)w+a_2(z)w^2+\cdots+a_{n}(z)w^{n}=0
\label{eqn:eqn001}
\end{equation}
with $z$ and $w$ complex variables and the coefficients, $a_i(z)$, polynomials in $z$ with rational coefficients.  By the Implicit Function Theorem, (\ref{eqn:eqn001}) defines locally, an analytic function $w(z)$ when $\displaystyle \frac{\partial f}{\partial w}\neq 0$.  The solution set of (\ref{eqn:eqn001}) defines an algebraic curve, $w(z)$, and it is known from the general theory of algebraic functions that $w(z)$ can be described in a disk centered at $z_0$ by $n$ fractional power series called Puiseux series with radii of convergence extending at least to the distance to the nearest singular point.  In an earlier paper by this author \cite{Milioto}, the method of Kaub and Traub was used to compute the power expansions, and a numeric integration method was described to compute their radii of convergence.  In this paper, a new method to compute radii of convergence is described, and the accuracy of the series is studied in their domain of convergence. 
 
A Puiseux expansion of (\ref{eqn:eqn001}) at a point $z_0$ is a set of $n$ fractional power expansion in $z$ given by 
\begin{equation}
 \{P_i(z)\}=\{w_i(z)\}=\sum_{k=r}^\infty a_k (z-z_0)^{\frac{m_k}{c}},\quad i=1,2,\cdots,n
\label{eqn:eqn002}
\end{equation} 
where $z$ lies in the domain of convergence of the series.  For finite $z_0$, $f$ is translated via $f(z+z_0,w)$ then $w(z)$ expanded as the set $\displaystyle \{w_i(z)\}=\sum_{k=r}^\infty a_k z^{\frac{m_k}{c}}$  with $z^{\frac{m_k}{c}}$  interpreted as principal-valued, and the series evaluated at the relative coordinate $z_r=z-z_0$.  

 In the case of an expansion at infinity, $f$ is translated via $\displaystyle g(z,w)=z^\delta f\left(\frac{1}{z},w\right)$ where $\delta$ is the largest exponent of $z$ in $f(z,w)$. Then an expansion of $w(z)$ at the origin in terms of $g(z,w)$ is an expansion of $w(z)$ at infinity with the series evaluated at $\displaystyle z_r=\frac{1}{z}$.   Section \ref{testCase4} is an example of an expansion at infinity.

The derivative of $w(z)$ at a point $(z,w)=(p,q)$ can be computed as
\begin{equation}
\frac{dw}{dz}=\lim_{(z,w)\to(p,q)} \left(-\frac{f_z(z,w)}{f_w(z,w)}\right)
\label{eqn:eqn003}
\end{equation}
when this limit exist. A singularity $\{z_s,w_s\}$  of $w(z)$ is a point where the limit does not exist and in this paper, the term \enquote{singular point} refers to the $z$-component of the singularity.
\section{Conventions used in this paper}
\begin{enumerate}
\vspace{1ex}\item  Computations are computed with a working precision of $1000$ digits.   An exception to this are the accuracy and order functions which do not need this level of precision and are only computed to machine precision.   All reported data however are shown with six or fewer digits for brevity.
\vspace{1ex}\item  Accuracy is the number of accurate digits to the right of a decimal point.  Precision is the total number of accurate digits in a number. For a number $x$,  $p=a+\log_{10}|x|$ with $p$ the precision and $a$ the accuracy.
\vspace{1ex}\item Finite singular points are the zeros of the resultant of $f$ with $\displaystyle \frac{\partial f}{\partial w}$ and are arranged into three lists:
\vspace{1ex}\begin{enumerate}
	\item \textbf{ The singular list:}  This is a list of the finite singular points in order of increasing distance from the origin.  Conjugate singular points are ordered real part first then imaginary part.  The singular points are labeled $s_1$ through $s_T$ with $T$ the total number in the list,
	\item \textbf{The singular sequence:}  For each singular point $s_b$, the remaining singular points are ordered in increasing distance from $s_b$,
	\item \textbf{The comparison sequence:}  This is a truncated singular sequence which is used to identify convergence-limiting singular points (CLSPs). 
\end{enumerate}

\vspace{1ex}\item Each singular point is assigned a circular perimeter with radius equal to 1/3 the distance to the  nearest singular point.   This circle is called the singular perimeter.
\vspace{1ex}\item A $k$-cycled branch refers to a part of $w(z)$, $k$-valued and  represented by $k$ Puiseux series with radius of convergence $R$.  
\vspace{10pt} \item The expansions of (\ref{eqn:eqn001}) at a singular point $s_b$ is a set of $n$ Puiseux expansions, $\{P_i(z)\}$ in terms of $\displaystyle z^{1/c}$ where $c$ is a positive integer and can be different for different series in the set.  $c$ is both the cycle size of the series and cycle size of the branch represented by the series.  For example, the power series $\displaystyle \sum_{k=r}^{\infty} a_k z^{\frac{m_k}{3}}$ has a cycle size of $3$ and represents a $3$-cycle branch:  The branch has three coverings over a deleted neighborhood of the expansion center and is represented by three such Puiseux expansions making a $3$-cycle conjugate set of power expansions.  The geometry of the branch is similar to the geometry of $\displaystyle z^{\pm 1/3}$. 

The set of $n$ expansions are numbered $1$ through $n$ with conjugate members sequentially numbered.  For example, a $15$-degree function may have a $5$-cycle branch with series numbers $\{7,8,9,10,11\}$, a $1$-cycle branch with series number $\{12\}$, a $3$-cycle branch with members $\{13,14,15\}$,  and a $6$-cycle with members $\{1,2,3,4,5,6\}$. 
\vspace{10pt} \item Branches of algebraic functions are categorized according to their algebraic and geometric morphologies into six types:  $T,E,F_p^q,V_p^q,P_p^q,L^q$. The $T,E$ and $L^q$ branches are $1$-cycle branches, and $F_p^q, V_p^q$ and $P_p^q$ are $p$-cycle branches.  Appendix \ref{appendix:appA} describes each branch type.
\vspace{10pt}\item Reference is made to a base singular point $s_b$.  This refers to a center of expansion of a Puiseux series with $s_b$ a singular point.  In the procedure described below, the branch surfaces about $s_b$ are  analytically continued over other singular points $s_n$ in order of increasing distance from the base singular point until the nearest convergence-limiting singular point (CLSP) is encountered.

Another closely related term is the impinging singular point or ISP of a branch sheet.  The ISP of a single-valued sheet of a multivalued branch is the nearest singular point impeding the analytic continuity of the branch surface.  The nearest ISP of all branch sheets in a conjugate set is the CLSP for the branch and establishes the radius of convergence of their power expansions.

The CLSP of a conjugate set is not unique as multiple (conjugate) singular points may impinge the analytic continuity of a branch sheet.  In these cases, the first member in the singular sequence is selected as the CLSP. 
\vspace{10pt}\item $R$ is a positive real number representing the radius of convergence of a power series centered at a singular point $s_b$.  The value of $R$ is expressed in terms of the associated CLSP.  For example, if a power expansion has a center at the tenth singular point $s_{10}$, and its CLSP was found to be $s_{25}$, then $R=|s_{10}-s_{25}|$.  This notation is presented as the exact symbolic expression for radius of convergence.
\vspace{10pt} \item The Puiseux expansions of $w(z)$ at a point $s_b$ are grouped into conjugate classes. For example, a $5$-cycle branch of $w(z)$ is expanded into five Puiseux series in powers of $z^{1/5}$, one series for each single-valued sheet of the branch.  These five series make up a single $5$-cycle conjugate class.  The sum of the conjugate classes at a point $z$ is always equal to the degree of the function in $w$.  A power expansion of a $10$-degree function consist of the set $\{P_i\}$ such that the sum of the conjugate  types is $10$.  This could consists of a single $10$-cycle conjugate class containing ten series, or three different $3$-cycle conjugate classes and a single $1$-cycle conjugate class or some other combination of conjugate classes adding up to $10$.  One member of each conjugate class  is selected as the class generator.  Each series member in a conjugate class can be generated by conjugation of a member of the class as follows:  

 Let
\begin{equation}
P_k(z)=\sum_{i=r}^{\infty} a_i z^{\frac{m_i}{c}}
\label{eqn:eqn300}
\end{equation}
be the $k$-th member of a $c$-cycle conjugate class of Puiseux series where all $\frac{m_i}{c}$ exponents are placed under a least common denominator $c$ and $\displaystyle z^{\frac{m_i}{c}}$ is the principal-valued root.  Then the $c$ members of this conjugate class are generated via conjugation of (\ref{eqn:eqn300}) as follows:
\begin{equation}
P_j(z)=\sum_{k=r}^{\infty} a_k \left(e^{\frac{2 j\pi i}{c}}\right)^{m_k} z^{\frac{m_k}{c}};\quad j=0,1,\cdots,c-1.
\label{eqn:eqn301}
\end{equation}
\vspace{10pt}\item  The order  of a $n$-cycle series of length $l$ is denoted by $\mathscr{O}$ and is the highest integer power of $z$ in the series nearest to the exponent of the $l$'th term.  Accuracy measurement in this study are done relative to a series order and not to a specific number of series terms and therefore accuracy results of multiple series will often include series of different lengths.  For example, $500$ terms of a $1$-cycle series can have an order of $500$ whereas $500$ terms of a $20$-cycle series may only attain an order of $30$ if the expansion has many fractional exponents in the series.  
\vspace{10pt}\item If an expansion of an $n$-degree function produces $n$ series in terms of $z^{1/n}$, the function fully-ramifies into a single $n$-cycle branch producing $n$ series belonging to an $n$-cycle conjugate class.  The branch is morphologically similar to $\displaystyle z^{\pm 1/n}$.  An $n$-degree function minimally-ramifies at a singular point if it ramifies into a single $2$-cycle branch and $(n-2)$ single-cycle branches. 
\vspace{10pt}\item  Absolute coordinates and relative coordinates in  the z-plane are used.    An absolute point $z_a$ is a  point in the z-plane.  A relative point $z_r$ is a point relative to a finite singular point $s_b$ given by
\begin{equation}
z_a=s_b+z_r,
\label{eqn:009}
\end{equation}
and in the case of an expansion at infinity, 
$$
z_a=\frac{1}{z_r}.
$$  
This is necessary for the following reasons:
\begin{enumerate}
	\item Power series are generated relative to an expansion center, $s_b$.  For example if the expansion center is $s_b=1+i$ and the series is evaluated at $z_r=0.25$, the accuracy is determined by computing a higher precision branch value $v_b$ by first solving for the roots $\{w_i\}$ of $f(1.25+i,w)=0$ and identifying which root $v_b$ corresponds to the series value at $z_r$.
\item The point $D$ in Figure \ref{figure:figure1} is an absolute point.  If $s_b=2+2i$ and $D=1.7-1.8i$, in order to evaluate a series expansions at $D$, the absolute point $D=1.7-1.8i$ is first converted to the relative point $z_r=1.7-1.8i-(2+2i)=-0.3-3.8i$. 
\item A series is evaluated over a list of points $C=\{z_i\}$ on a circle around the singular point $s_b$ to compute accuracy profiles.  The points $C=\{z_i\}$ are relative coordinates to $s_b$.  The points in $C$ are translated to absolute coordinates as $C_a=\{s_b+z_i\}$, then the roots $\{w_i\}$ of $f(w,C_a)=0$ are computed to a higher precision and corresponding branch values identified to determine series accuracies.
\item An expansion at infinity is generated relative to an expansion at zero, and the associated power expansions use relative coordinates $\displaystyle z_r=\frac{1}{z_a}$.  See Section \ref{testCase4}.
\end{enumerate}
\vspace{10pt}\item The series comparison test used to determine CLSPs relies on comparing a base series value $v_s$ of a branch expansion at a point $D$ in Figure \ref{figure:figure1} to a list of series values $\{u_i\}$ computed at the next nearest singular point $s_n$ at point $D$.  Since this involves comparing numeric values at finite precision, a separation tolerance $s_t$ is used to identify a match.  $s_t$ is $1/10$'th the minimum separation of the members in $\{u_i\}$.  If $|v_s-u_k|<s_t$, then the $k$'th series value at $s_n$ is a match for $v_s$.  If this tolerance is exceeded as in the case of branch values being very close to on another, or an insufficient number of terms in the series, or multiple matches, the series comparison test halts and the analysis reverts to the numerical integration method.  This is described further in Section \ref{sec:sec006}
\vspace{10pt}\item An accuracy profile of a conjugate set of series is generated by computing the accuracy of generator series over a region in their domain of convergence.  For each value $|z|< R$, the accuracy is determined by comparing the series results to more precise roots $\{w_i\}$ of $f(z,w)=0$.  These roots are computed with Mathematica's \texttt{NSolve} function. This accuracy is then fitted to an accuracy function $A(r_f,o)$.
\vspace{10pt}\item Table \ref{table:table23} is a list of symbols used in this paper.

\begin{table}[ht]
\caption{Symbols}
$
\begin{array}{|c|l|}
\hline
\text{Symbol} & \text{Description}\\
\hline
m & \text{Number of series terms used in a calculation}\\
M & \text{Maximum number of terms in a series} \\
R & \text{Radius of convergence} \\
r_f & \text{A positive rational number between 0 and 1}\\
r & \text{Radius of $z=re^{it}$}\\
s_a & \text{Series accuracy}\\
A(r_f,o) & \text{Accuracy function. (Section \ref{section:section005})} \\
e_a & \text{Expected accuracy at $|z|< R$ given by $A(r_f,o)$}\\
\mathscr{O}(r_f,e_a) & \text{Order function (Section \ref{section:section92})}\\
(a,b,c,d)& \text{Coefficients of accuracy profile function}\\
s_b & \text{Singular point at center of expansion} \\
s_n & \text{Singular point $n$ in the singular list} \\
z_a & \text{Absolute coordinate value of $z$} \\
z_r & \text{Relative coordinate of $z_a$ with respect to a singular point} \\
\{w_i\} & \text{$n$ values of $w(z_a)$ computed to high precision} \\
\{v_b\} & \text{Particular set of branch values with  $\{v_b\}\subseteq \{w_i\}$}\\
\{v_s\} & \text{Branch series values at $z_r$ corresponding to the set $\{v_b\}$}\\
c_e & \text{Comparison error given by $c_e=|w_i-v_s|$}\\
r_t & \text{Residual tolerance. (Section \ref{section001})}\\
s_t & \text{Separation tolerance.  (Section \ref{sec:sec006})}\\
m_s & \text{Minimum separation. (Section \ref{sec:sec006})}\\
s_f & \text{Separation factor.  (Section \ref{sec:sec006})}\\
c_z & \text{Coefficient zero.  (Section \ref{sec:sec009})}\\
\{c_i\} & \text{Roots to the polygonal characteristic equation}\\
N_{zm} & \text{Maximum number of zero modular terms. (Section \ref{section:section109}) }\\
\mathscr{G} & \text{Riemann surface genus. (Section \ref{sec:sec500})}\\
\mathscr{K} & \text{Riemann-Hurwitz sum (Section \ref{sec:sec500})}\\
\hline
\end{array}
$
\label{table:table23}
\end{table}

For example, let $z=r_f R e^{it}$ with $R$ the radius of convergence of a series.  $v_s$ is the value of the series at $z$.  $v_b$ is the value of the corresponding branch at $z$ computed to a  higher precision than the series precision.  $c_e=|v_b-v_s|$ is the comparison error.  The accuracy of the series then becomes the negative of the exponent of $c_e$. 

Generator series are evaluated in their domain of convergence along circular domains $|z|< R$ and the accuracy determined by comparing to the corresponding member in the set $\{w_i\}$. The accuracies are then fitted to an accuracy function $A(r_f,o)$ which gives the expected accuracy, $e_a$, of the series as a function of the radial ratio $r_f$ and order $o$ of the series.  Solving for $o$ in $e_a=A(r_f,o)$ gives $o=O(r_f,e_a)$ which is the order function for estimating the order of a series needed for a desired accuracy $e_a$ at $z=r_f R e^{it}$.

The genus, $\mathscr{G}$, of $w(z)$ is easily calculated via the Riemann-Hurwitz formula once conjugate classes at all singular points are found.  The Riemann-Hurwitz sum $\mathscr{K}$ must be an even number and serves as a necessary (but not sufficient) check-sum of the overall cycle geometry.  See Section \ref{sec:sec500}.  
 \vspace{10pt}\item Power expansions are generated by the Newton Polygon method \cite{Kung}.  This algorithm has two types of function iterations:
\begin{enumerate}
	\item \textbf{Polygon iteration: } The first step in the algorithm is to create a Newton polygon establishing the initial terms of each expansion.  If there are multiple roots in the resulting characteristic equation, a new Newton polygon is created by iteration via the expression
	\begin{equation}
	f_{2}=z^{-\beta_i}f\left(z,z^{\lambda_i}(w+c)\right)
	\label{eqn:eqnpi1}
	\end{equation}
	and a second Newton polygon created for $f_{2}$.  If the resulting characteristic equation has multiple roots, polygon iterate $f_3$ of $f_2$ created and so on until the characteristic roots are simple.  The substitution $f\left(z,z^{\lambda_i}(w+c)\right)$ can cause numerical errors if not pre-processed beforehand.  See Section \ref{sec:sec009}. 
	\item \textbf{Newton Iteration: } Upon obtaining simple characteristic roots, the final polygon iterate $f_n$, after two additional transformations, is iterated by a Newton-like iteration to produce the desired number of expansion terms via the expression
	\begin{equation}
	   w_{j+1}=w_j-\text{mod}\left(\frac{\overline{f}(z,w_j)}{\overline{f}_w(z,w_j)},z^{2^{(j+1)}}\right);\quad w_0=c_i.
  	\end{equation}
		
		In the case of fractional polynomial solutions, the modular function continues to return zero after reaching the polynomial.  Finite polynomial solutions however have to be distinguished from infinite solutions with extremely large gaps in exponents between successive series terms which would also return zero for a (often small) number of iterations.  This is done by setting the number of maximum modular zeros, $N_{zm}$ to a large but manageable number.  In this case, $N_{zm}$ was set to $15$. 
		
		See \href{https://arxiv.org/pdf/2111.11883.pdf}{Determining radii of convergence of fractional power expansions around singular points of algebraic functions} for additional information about these concepts.
		\end{enumerate}
		
\end{enumerate}
  This study is organized as follows:
\begin{enumerate}
  \item Compute the singular points to $1000$ digits of precision,
	\item Compute conjugate classes at all singular points,
	\item Compute at least $1000$ terms of each branch expansion at a base singular point $s_b$ with a working precision sufficient to obtain series with at least $900$ digits of precision,
	\item Compute expansions for all singular points in the comparison sequence to at least $100$ terms,
	\item Estimate radius of convergence for each branch around $s_b$ using the Root Test,
	\item Compute CLSPs via the comparison test and integration test,
	\item Compute the accuracy function $A(r_f,o)$ and order function $O(r_f,e_a)$,
	\item Generate convergence, accuracy, and timing data for six test functions.
\end{enumerate}
\section{Precision of computations}
\label{sec:sec009}
The precision of a series is limited by the precision of the singular points as well as the reduction in precision incurred by various steps in the Newton polygon procedure. 
\begin{figure}[ht]
\centering
\includegraphics[scale=0.75]{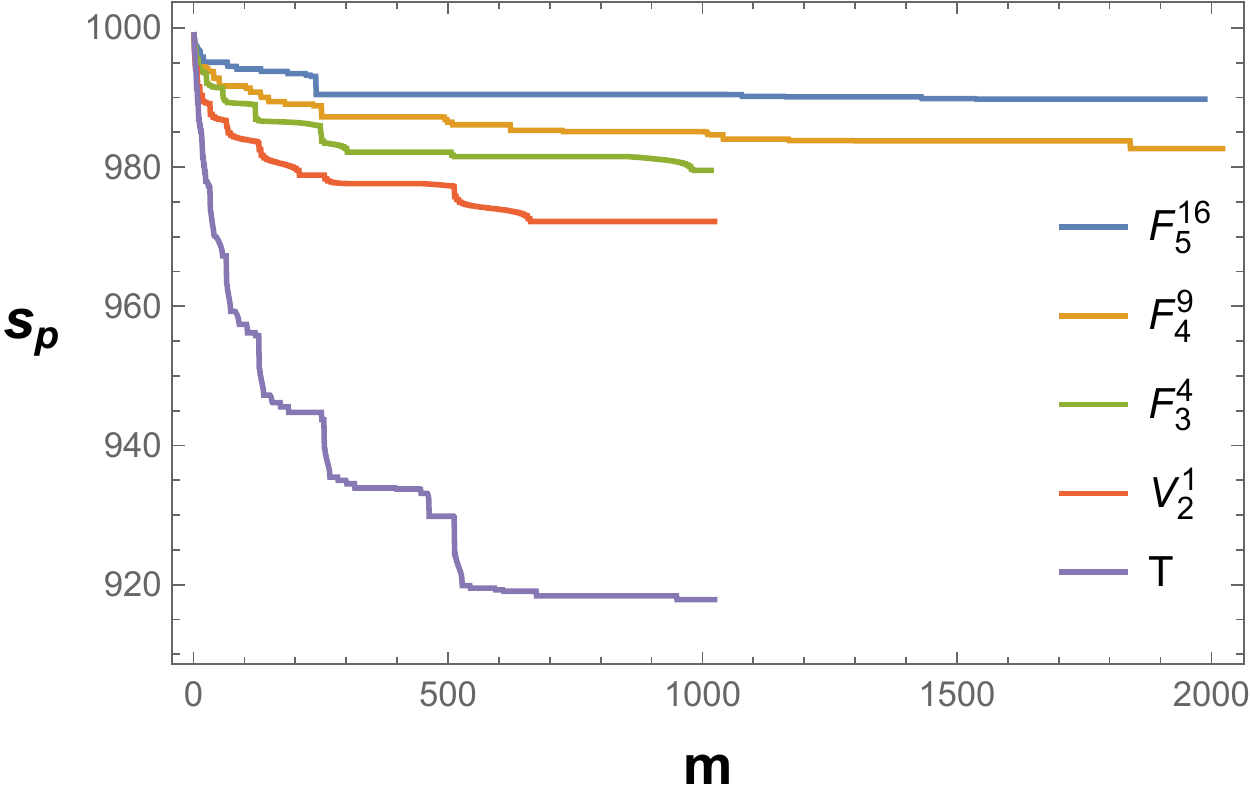}
\caption{Precision Profile Plot}
\label{figure:figure14}
\end{figure}

 Figure \ref{figure:figure14} is a precision plot of the five generator series for Test Case $1$ and exhibits a reduction in precision after each Newton iteration.  For example, the $1$-cycle $T$ branch shown as the purple line begins with $1000$ digits of precision  dropping to $918$ digits at the $1024$'th term.  The precision of each series is therefore dependent on the the number of terms used.  If the $1T$ series was evaluated very close to its center of expansion, the accuracy could rapidly increase as additional terms are added to the computation.  At some point, the precision of the results could approach the maximum precision of the series and no longer increase as more terms are used  leading to a case in which the  precision of the comparison error $c_e=|w_i-v_s|$ drops to zero.  This would skew the accuracy results and subsequently the fit functions $A(r_f,o)$ and $O(r_f,e_a)$.  To avoid this possibility, $c_e$ values with zero precision are omitted from the accuracy results.

 The Newton polygon algorithm produces translated functions $f(z+s_i,w)$ and $f(z,w+c_i)$.  These substitutions can create coefficients which are actually zero but due to finite numerical precisions, result in very small residual values.  Two residual conditions arise which are processed in the following order: 
 
\begin{enumerate}
	\item \textbf{R1:}  A singular point can be a zero to a coefficient of the transformed function $f(z+s_i,w)$.   These are called coefficient zeros denoted by the symbol $c_z$.  Since the singular points are computed to a finite number of digits, coefficients which are actually zero can have very small non-zero residues.  If not eliminated, this would cause incorrect polygon calculations.  The two sets $\{p_i\}$ and $\{q_i\}$, described below, are such singular points.  There may however exist other singular point zeros.  $f(z,w)$ is pre-processed to remove these zero coefficient before the substitution $f(z+s_i,w)$.  

	\item \textbf{R2:} Roots $\{c_i\}$, of the polygonal characteristic equations are substituted into the polygonal iterates as $f\left(z,z^{\lambda}(w+c_i)\right)$ and likewise results in zero coefficients which may have very small residues and must be removed. Similar to the \textbf{R1} pre-processing, these coefficients are first identified and removed prior to the substitution $w\to w+c_i$.

\end{enumerate}
\subsection{Additional precision issues}
\begin{enumerate}
	\item  Singular points with very small absolute values can lead to very small coefficients in the translated function $f(z+s_i,w)$ if the function has high powers of $z$.  For example if $s_i=1/1000$, then the substitution $z\to z+s_i$ into $z^{50}$ leads to a coefficient on the order of $10^{-150}$, and this small value must not be lost due to inadequate numerical precision as doing so would adversely affect the polygonal iteration step of the Newton polygon algorithm producing incorrect initial segments.  Since the test cases below were run with an average precision of $1000$, this particular cases would be correctly processed.  However, there exists functions with arbitrarily small singular sizes which would be mis-handled.  The singular size therefore is carefully monitored and if the size exceeds the precision limitation of the calculation, the procedure is halted. 
	\item \textbf{Multiple polygon iteration}  Each polygon iterates of (\ref{eqn:eqnpi1}) can reduce the precision of the resulting function iterate $f_i$.  
			\item \textbf{Limitations of Mathematica's \texttt{NSolve} function: }Roots returned by \texttt{NSolve} are limited in precision to the precision of the input equations.  If the equations are at $1000$ digits of precision, then the maximum precision of the roots is $1000$.  However, the precision of the roots returned by \texttt{NSolve} may be less than the precision of the equations.  For example, the roots of $1/2+x^2$ returned by \texttt{NSolve} when the precision of the equation is set to $1000$ is $1000$.  However, if the precision of $1/4+x+x^2$ is set to $1000$, \texttt{NSolve} returns roots to only $500$ digits of precision.  
				\end{enumerate}

\section{Computing the singular points}
\label{singularPointsSection}
The finite singular points are computed by solving for the zeros of the resultant of $f$ with $\displaystyle\frac{\partial f}{\partial w}$  using Mathematica's \texttt{NSolve} function.  This computation is CPU-intensive when $f$ is of high degree and the coefficients $c_i(z)$ are non-sparse and high degree.  The $50$-degree function studied in Test Case $4$ took $3.8$ hours to compute $4584$ singular points to $1000$ digits of precision, whereas a random $20$-degree function with $156$ singular points and low degree coefficients takes about one second. 
   
 There are two sets of singular points by inspection:
\begin{enumerate}
	\item Roots of $a_0(z)$ if $a_1(z)=0$:  These are the set $\{q_i\}$,
		\item Roots of $a_n(z)$:  These are the set of poles $\{p_i\}$.
		\end{enumerate}
\section{Computing initial terms and identifying conjugate class membership}
\label{section001}
The method of Kung and Traub \cite{Kung} implements Netwon polygons to generate the initial terms (segments) of each power expansion.  See also \cite{Milioto} for more information about Newton polygons.  For a generic polynomial, computing the initial terms requires only a few iterations of Newton polygon and is executed quickly.  Even  the $50$-degree polynomial in Test Case $5$, required at most two seconds to compute the initial segments at a singular point. 

 However, in many cases, not all initial segments require further expanding.   Rather, in this paper the initial terms are first used to identify the cycle size of each expansion and their conjugate class membership.  Class membership is obvious from the initial segment exponents when there are distinct multi-cycles.  In the case of multiple $k$-cycle segments, membership is determined by conjugating the initial terms in order to determine which expansions belong to each conjugate set.  One member of each class is selected as the class generator and further expanded via Kung and Traub's method of iteration. 
 
Consider an expansion at the origin of the following $12$-degree function:
\begin{equation}
f(z,w)=w^4 z \left(\text{a2} \left(1-w^2 \left(\text{a2}^*\right)^2\right)\right)^4-\left(\text{a2}^* \left(\text{a2}^2-w^2\right)\right)^4;\quad a=2/3+i/4,
\end{equation}
and the initial terms obtained from the Newton Polygon step given by $(\ref{eqn:eqn2006})$.  Since there are $12$ expansions all of which are $4$-cycle, there are three $4$-cycle conjugate classes of branch expansions:  
\begin{equation}
\begin{aligned}
P_1(z)&=-0.667-0.25i-(0.28+0.244i) z^{1/4} \\
P_2(z)&=-0.667-0.25i-(0.244-0.28i) z^{1/4}\\
P_3(z)&= -0.667-0.25i+(0.244-0.28i) z^{1/4}\\
P_4(z)&=-0.667-0.25i+(0.28+0.244i) z^{1/4}\\
P_5(z)&= 0.667+0.25i-(0.28+0.244i) z^{1/4} \\
P_6(z)&=0.667+0.25i-(0.244-0.28i) z^{1/4}\\
P_7(z)&= 0.667+0.25i+(0.244=0.28i) z^{1/4} \\
P_8(z)&= 0.667+0.25i+(0.28+0.244i) z^{1/4}\\
P_9(z)&=\displaystyle -\frac{1.973}{z^{1/4}} \vspace{8pt}\\
P_{10}(z)&=\displaystyle  -\frac{1.973i}{z^{1/4}} \vspace{8pt} \\
P_{11}(z)&= \displaystyle \frac{1.973 i}{z^{1/4}} \vspace{8pt} \\
P_{12}(z)&= \displaystyle \frac{1.973}{z^{1/4}}.\\
\end{aligned}
\label{eqn:eqn2006}
\end{equation}
In order to determine which expansions belong to each 4-cycle conjugate class, segment members are conjugated.  Consider:
$$
P_1(z)=(-0.666667-0.25 I)-(0.2799 +0.244276 I) z^{1/4}.
$$
Conjugation of $P_1$ produces the following list of members:
\begin{equation}
\begin{array}{c}
(-0.666667-0.25 i)-(0.2799 +0.244276 i) z^{1/4}\\
(-0.666667-0.25 i)+(0.244276 -0.2799 I) z^{1/4}\\
(-0.666667-0.25 i)+(0.2799 +0.244276 I) z^{1/4}\\
(-0.666667-0.25 i)-(0.244276 -0.2799 I) z^{1/4},\\
\end{array}
\end{equation}
and these are $P_1$ through $P_4$.  Therefore, $P_1$ through $P_4$ are the four members of a $4$-cycle conjugate class $1V_4^1$ and the series numbers for this set are $\{1,2,3,4\}$.  Conjugating expression $P_5(z)$ in \ref{eqn:eqn2006} gives the next four series, $\{5,6,7,8\}$ making up a second $4$-cycle class $2V_4^1$, and conjugating $P_9$ gives the set of the next four series, $\{9,10,11,12\}$ making up a third $4$-cycle conjugate set $3P_4^{-1}$.  Series $1$, $5$, and $9$ are selected as the generators of the three conjugate classes and further expanded via Newton Iteration.  Once the desired number of terms for each generator series has been computed, the full $12$ expansions around this singular point can be generated by conjugating the generator series.  Since conjugation is much faster than Newton iteration, computing the $12$ series this way is faster than generating each member separately via Newton Iteration.

Computing $1024$ terms of each generator series at $1000$ digits of precision took $3.4$ minutes.  Conjugation of all generators took one second.  Compared to expanding all initial segments, this represents a four-fold reduction in execution time.  The performance gain is dependent upon the ramification extent at the expansion center with minimal gain obtained with minimal ramification.

\section{Expanding the initial segments via Newton-like iteration}
\label{section:section109}
The following is a brief summary of the Kung and Traub method of iterating the initial series segments.  For more information, interested readers are referred to the author's website: \href{https://www.jujusdiaries.com/p/section-5.html}{Examples of power expansions around singular points of algebraic functions.} which includes worked examples. 

Consider the following function from the website:
\begin{equation}
f(z,w)=(z+z^2)+(1+z)w+w^2=0.
\end{equation}
Processing this function through the Newton Polygons twice produces the first two terms (initial segments) of each series 
\begin{align*}
P_1(z)&=-2/3-i z^{1/2} \\
P_2(z)&=-2/3+i z^{1/2},
\end{align*}
and the second function iterate $f_2(z,w)$. 

 A critical part of a numerical Newton polygon algorithm is accurately identifying multiple roots of the characteristics equation. This is accomplished by first setting the roots to the same precision.  When two numbers accurate to a set precision are subtracted in Mathematica, the precision of the difference is zero.  Thus, the roots are first set to the minimum precision of the set, and the differences between roots are checked with those having zero precision identified as multiples.  Although the roots are computed with a default precision of $1000$ digits, the actual precision can be significantly lower as described above.  In these cases, the singular points are generated with a sufficient precision to obtain roots near $1000$ digits of precision.

As the Newton polygon algorithm will often lead to fractional polynomials, $f_2$ is transformed into a polynomial with integer powers by the following two transformations:
\begin{equation}
\begin{aligned}
\widehat{f}(z,w)&=\frac{1}{z^{\beta_k}} f_k(z,z^{\lambda_k} w) \\
\overline{f}(z,w)&=\widehat{f}(z^d,w)
\end{aligned}
\end{equation}
where $k$ is the index of the last polygonal function which did not produce a characteristic equation with multiple roots ($k$ is $2$ in this case). Let $d$ be the lowest common denominator of the exponents $\{\lambda_j\}$ for $j=1,2,\cdots,k$ for each of the segments above.  In this case, $k=2$ and $d=2$ with $\lambda_2=1/2$ and $\beta_2=1$.  Therefore we have:
\begin{equation}  
\begin{aligned}
\widehat{f}(z,w)&=\frac{1}{z^{\beta_k}} f_k(z,z^{\lambda_k}w)  \\
   &=z^{-1}\left[(z+z^2)+z(z^{1/2}w)+zw^2\right] \\
   &=1+z+z^{1/2}w+w^2=0
\end{aligned}
\end{equation}
and
$$
\begin{aligned}
\overline{f}(z,w)&=\widehat{f}(z^d,w)\\
   &=1+z^2+zw+w^2.
\end{aligned}
$$
    In order to generate more terms of the series, Kung and Traub implements a Newton-like iteration on $\overline{f}$:
   \begin{equation}
   w_{j+1}=w_j-\text{mod}\left(\frac{\overline{f}(z,w_j)}{\overline{f}_w(z,w_j)},z^{2^{(j+1)}}\right);\quad w_0=\{i,-i\}
	\label{eqn:eqn989}
   \end{equation}
where the $\text{mod}$ function extracts all terms of the Taylor expansion of the quotient with power less than $2^{(j+1)}$.  Listing $1$ is the modulus step in (\ref{eqn:eqn989}) implemented in Mathematica.
\begin{center}
\begin{minipage}{0.85\linewidth}
\begin{lstlisting}[style=myMathematica,
frame=single,
caption=Mathematica code,
label=code4]
newtonIterate = 
  Normal[Series[fBar/fBarDeriv, {z, 0, 2^(j + 1) - 1}]];
	\end{lstlisting}
	\end{minipage}
	\end{center}
	
	Solutions with fractional polynomial solutions return sequential zero modular result after a finite number of iterations.  In this study, this number of maximum zero modular values is $N_{zm}$ and set to $15$.  See Test Case \ref{testCase2} for an example of a polynomial solution.
%
%
\section{Approximating radii of convergence via the Root Test}
  The Root Test is used to approximate the radius of convergence of a branch expansion in order to estimate the number of additional power expansions needed to identify the CLSP of the branch series.  The Root Test however is not applicable to polynomial solutions which have infinite radii of convergence.  

 The standard definition is modified to include the branch cycle size:
\begin{equation}
R=\frac{1}{\displaystyle \liminf_{k\to \infty}\; |a_k|^{\frac{c}{m_k}}}
\label{eqn:eqnss2}
\end{equation}
where $c$ is the cycle size of the series, and the set $\{m_i\}$ is the set of exponent numerators under a least common denominator.  For example, the terms $a_0+a_1 z+a_2 z^{3/2}+a_3 z^{9/4}+z^{3}$ would have the set $\{m_i\}=\{0,4,6,9,12\}$.  Then the radius of convergence of each branch expansion can be approximated by forming the set:
$$
S=\bigg\{\left(\frac{1}{m_k},\frac{1}{|a_k|^{\frac{c}{m_k}}}\right)\bigg\}
$$
and extrapolating $\displaystyle \lim_{k\to\infty} S$ using a sufficient number of trailing points of $S$. 
For example, consider $512$ terms of a $2$-cycle series:
$$
a_1 z^{1/2}+a_2 z^{2}+a_3 z^{5/2}+a_4 z^{3}+\cdots+a_{512}z^{512}.
$$
Therefore
$$
S=\biggr\{\left(\frac{1}{1},\frac{1}{|a_1|^{1/2}}\right),\left(\frac{1}{4},\frac{1}{|a_2|^{1/2}}\right),\left(\frac{1}{5},\frac{1}{|a_3|^{2/5}}\right),\left(\frac{1}{6},\frac{1}{|a_3|^{1/3}}\right),\cdots,\left(\frac{1}{1024},\frac{1}{|a_{512}|^{1/512}}\right)\biggr\}.
$$
The radius of convergence, $R$, is approximated by minimizing a suitable curve to the greatest lower bound of $S$ using Mathematica's  \texttt{Minimize} function and then extrapolating to zero.  Linear, quadratic and cubic curves were used as test fit functions with the smallest residual error selected as the best fit.  Figure \ref{figure:plot1} shows the generator series points for the $F_5^{16}$ branch of Test Case $1$ as blue points.  The dashed red curve is the best fit of the greatest lower bound of points.  The black point is the extrapolated value of $0.645$ for the radius of convergence.  The actual radius of convergence for this branch is $|s_{27}|\approx 0.641$.  The  CLSP for the power expansion are approximated by selecting a singular point with absolute value closest to the extrapolated point and then used to estimate the minimum number of singular points in the comparison sequence. 
\begin{figure}[ht] 
\centering
\includegraphics[scale=0.75]{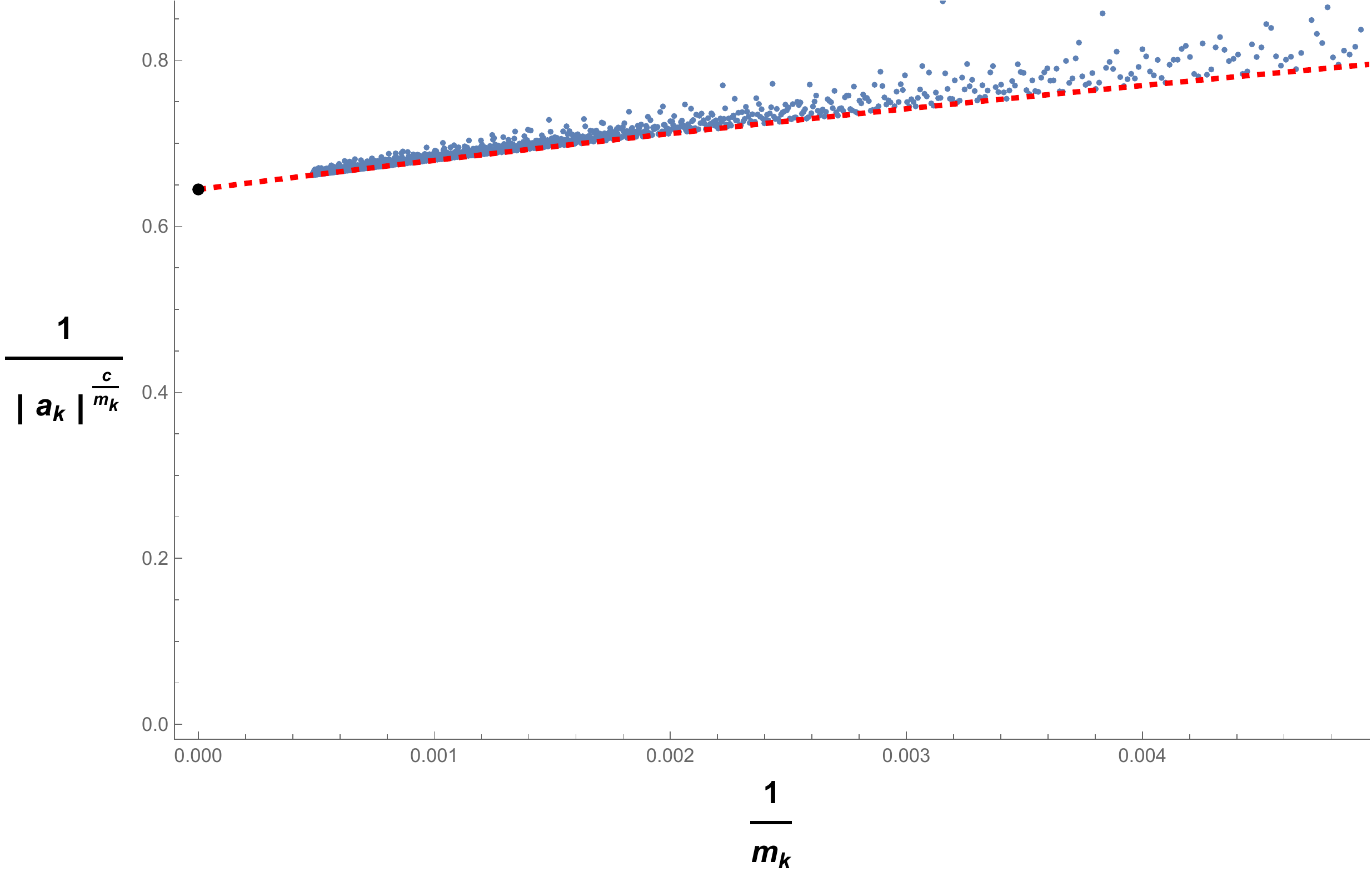}
\caption{Root Test results of $F_5^{16}$ branch of Test Case 1}
\label{figure:plot1}
\end{figure}
\section{Computing the comparison sequence}
 The most distant estimated CLSP determined by the Root Test determines the minimum length of the comparison sequence used by the series comparison and integration tests to identify CLSPs at a base singular point $s_b$.  Table \ref{table:table19} gives the Root Test results for Test Case $1$, and since the expansion center is the origin, $s_{118}$ is the most distant estimated CLSP giving a minimum comparison sequence $s_2$ through $s_{118}$.
 
Expansions must be generated at all singular points in the comparison sequence but do not require a large number of terms since these series are evaluated at their singular perimeter with good convergence.  Since the Root Test is only an approximation for the most distant CLSP, a reasonable comparison sequence in this case is $s_2$ through $s_{125}$.

\begin{table}[ht]
\caption{Test Case $1$ Root Test results at the origin}
$
\begin{array}{|c|c|c|c|c|c|c|}
\hline
 \text{Index} & \text{Sheet} & \text{Type} & \text{Cycle Size} & \text{est. R} & \text{est. CLSP} & \text{Series size} \\
\hline
 1 & 1 & \, F_5^{16} & 5 & 0.64571 & 31 & 1987 \\
 2 & 6 & \, F_4^9 & 4 & 0.50751 & 7 & 2022 \\
 3 & 10 & \, F_3^4 & 3 & 0.16797 & 2 & 1017 \\
 4 & 13 & \, V_2^1 & 2 & 0.16772 & 2 & 1024 \\
 5 & 15 & \text{T} & 1 & 1.0966 & 118 & 1024 \\
\hline
\end{array}
$
\label{table:table19}
\end{table}

\section{Computation of convergence-limiting singular points}
\subsection{Necessary and sufficient conditions for analytically-continuing a branch across a singular point:}
\label{sub:section0041}
\begin{enumerate}
	\item  In order to analytically continue a $k$-cycle branch from one singular point to the next nearest singular point, the next singular point must have at least $k$ single-cycle analytic branches to support continuity, i.e., $1$-cycle branches which do not have poles,
	\item A $k$-cycle branch is continuous across a singular point if all branch sheets continue onto analytic $1$-cycle branches, 
	\label{enum:enum1}
\end{enumerate}

Individual single-valued branch sheets of an $k$-cycle branch may continue across different singular points but the analytic region of the branch as well as the convergence domain of its power expansion is established upon analytic continuity with the nearest multi-cycle branch sheet or branch sheet with a pole.  The first singular point in which this occurs is the CLSP for the associated set of conjugate Puiseux series and establishes their radii of convergence. 
\section{Methods to compute CLSPs}
 Two methods are used to find CLSPs:
\subsection{Constructing an analytically-continuous route between singular points by numerical integration}
\label{sec:sec005}
 Since the Puiseux series for (1) converge at least up to the next nearest singular point, an analytically-continuous route can be created between the perimeter of the base singular point and perimeter of the next nearest singular point.  This is shown in Figure \ref{figure:figure1}.  In the diagram the circles are the singular perimeters.  A straight line path between point $A$ and $D$ is given by the expression
\begin{equation}
z(t)=A(1-t)+D t;\quad 0\leq t\leq 1,
\label{eqn:eqn154}
\end{equation}
and then each branch sheet of $s_b$ is  numerically integrated over the path from point $A$ to point $D$ onto a branch sheet of $s_n$via the following set of initial value problems:
\begin{equation}
\frac{dw}{dt}=-\frac{f_z}{f_w}\frac{dz}{dt}, \quad\; w_i(0)=P_i(A-s_b); \quad i=1,2,\cdots,n 
\label{eqn:eqn254}
\end{equation}
with  $z(t)$ defined by $(\ref{eqn:eqn154}), 0\leq t\leq 1$, and each $P_i(z)$ is a Puiseux series at $s_b$.  Each $s_b$ branch sheet is then checked for analytic continuity over $s_n$ as per Section \ref{sub:section0041}.  See also \cite{Milioto3} for further background about this technique.

 After $s_n$ has been checked and branch sheets have been found to be continuous, a path to the next sequential singular point is created and the previously-continued branches tested for further continuity.  However there is the possibility of attempting to continue a branch sheet to another singular point when a removable singular point is in the path of integration.  Numerical integration will fail over a removable singularity even though the function is analytic because Equation (\ref{eqn:eqn002}) is used for the derivative and at a removable singular point, this quotient is indeterminate.   In this case, the integration path is split into paths $\beta_1$, around half the perimeter of the removable singular point $E$, and over $\beta_2$ shown in the second diagram of Figure \ref{figure:figure1}.
\begin{figure}[ht]
\centering
\includegraphics[scale=0.75]{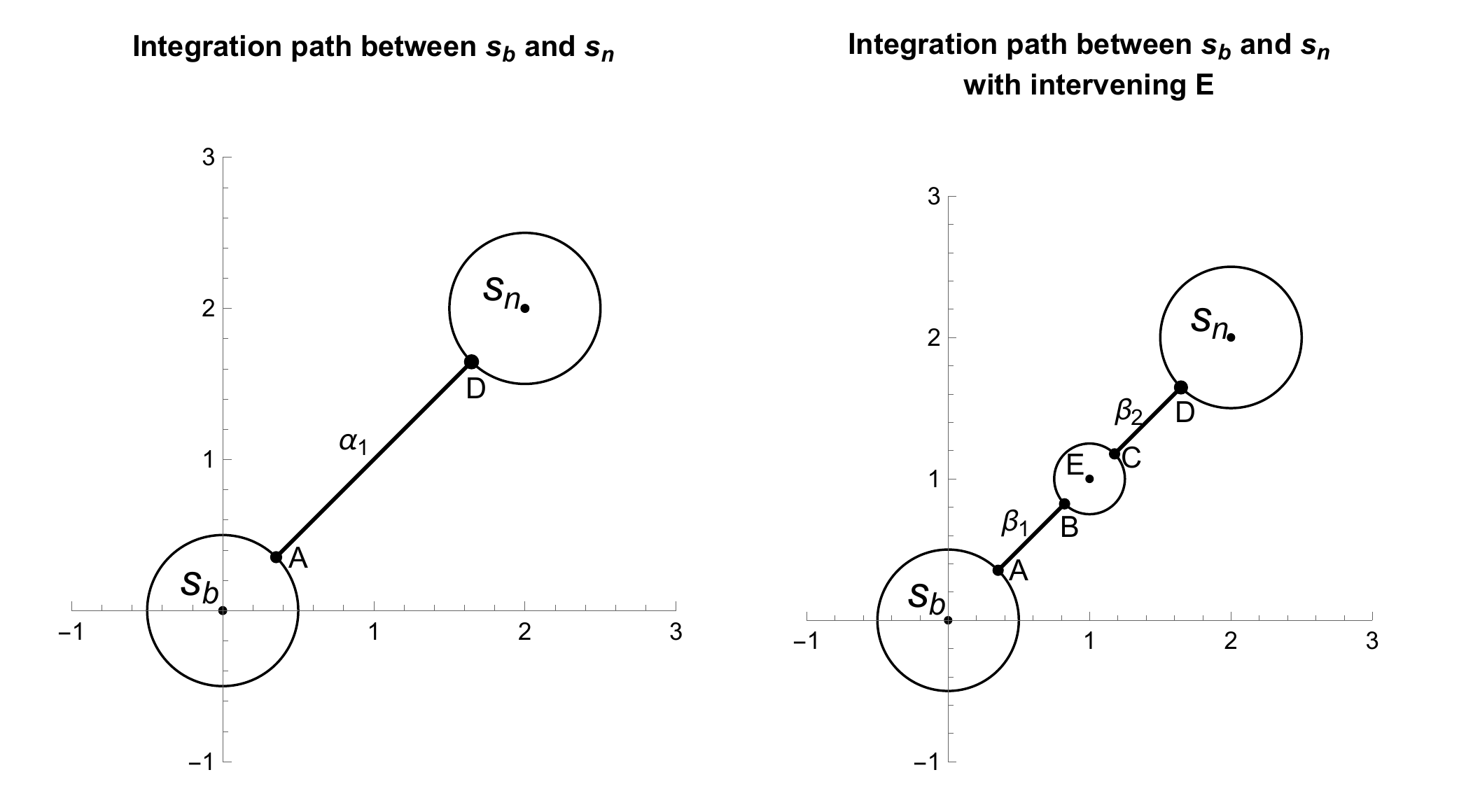}
\caption{Continuation path used by numerical integration method}
\label{figure:figure1}
\end{figure}
\subsection{Identifying CLSPs by the series comparison test}
\label{sec:sec006}

This section introduces a simpler method for computing a CLSP based on a series comparison test. 

 In the series comparison test,  the value of a base branch sheet at point $D$ in Figure  \ref{figure:figure1} is compared  to the expansions centered at $s_n$ at point $D$ using the separation threshold $s_t$.  Convergence of both series is guaranteed since the series centered at $s_n$ converges on its perimeter, and convergence of the current analytically-continuous base series are guaranteed to converge at $D$ since they converge at the previous $n-1$ singular points.  However, the convergence rate of the base sheets will decrease as successive singular points approach the CLSP of a base expansion sheet.  This necessitates computing the base series at a sufficiently high precision and number of terms. 
\begin{table}[ht]
	\caption{Branch sheet values at $s_2$}
	$
	\begin{array}{|l|}
	\hline
 -2.44964-1.22342 i \\
 -1.0006+1.0348 i \\
 -0.92827-0.314806 i \\
 -0.783318+0.21872 i \\
 -0.696357-0.639337 i \\
 -0.371485-0.178941 i \\
 0.371485\, +0.178941 i \\
 0.696357\, +0.639337 i \\
 0.783318\, -0.21872 i \\
 0.92827\, +0.314806 i \\
 1.0006\, -1.0348 i \\
 2.44964\, +1.22342 i \\
\hline
\end{array}
$
\label{table:table50}
\end{table}

Using the data from Section \ref{section001}, with three $4$-cycle branches and three singular points, consider series $1$ expansion at the origin with a value of $v_s=-0.92827-0.314806 i$ at point D in Figure \ref{figure:figure1} and compare it with the $12$ expansion values at $s_2$ at point D given in Table \ref{table:table50}.  In this case, there are no poles at $s_1$ and $s_2$.   The minimum distance between the points in Table \ref{table:table50} is $0.04$ so that the separation threshold is $s_t=0.04/10=0.004$ and clearly $v_s$ is within this threshold for the third entry in this list.  The actual difference between the two values in this test was $10^{-40}$.  Therefore, the branch sheet corresponding to series $1$  continues onto the branch sheet associated with the third series at $s_2$. 

		As the analysis moves farther away from the base expansion and nearer to the CLSP, the accuracy of the base expansions will decrease and may not satisfy the separation threshold or may produce multiple matches.  This condition is checked in the algorithms and when encountered, the comparison test is halted and flagged for the numerical integration test. 
\section{Illustrative example of Radius of convergence computation using series comparison test}
\label{section:section92}
Consider the function from Test Case $1$:
\begin{equation}
f_1(z,w)=\left(z^{30}+z^{32}\right)+\left(z^{14}+z^{20}\right) w^5+\left(z^5+z^9\right) w^9+\left(z+z^3\right) w^{12}+(6)
w^{14}+\left(2+z^2\right) w^{15}.
\label{eqn:eqn1009}
\end{equation}

\begin{figure}[ht] 
\centering
\includegraphics[scale=0.75]{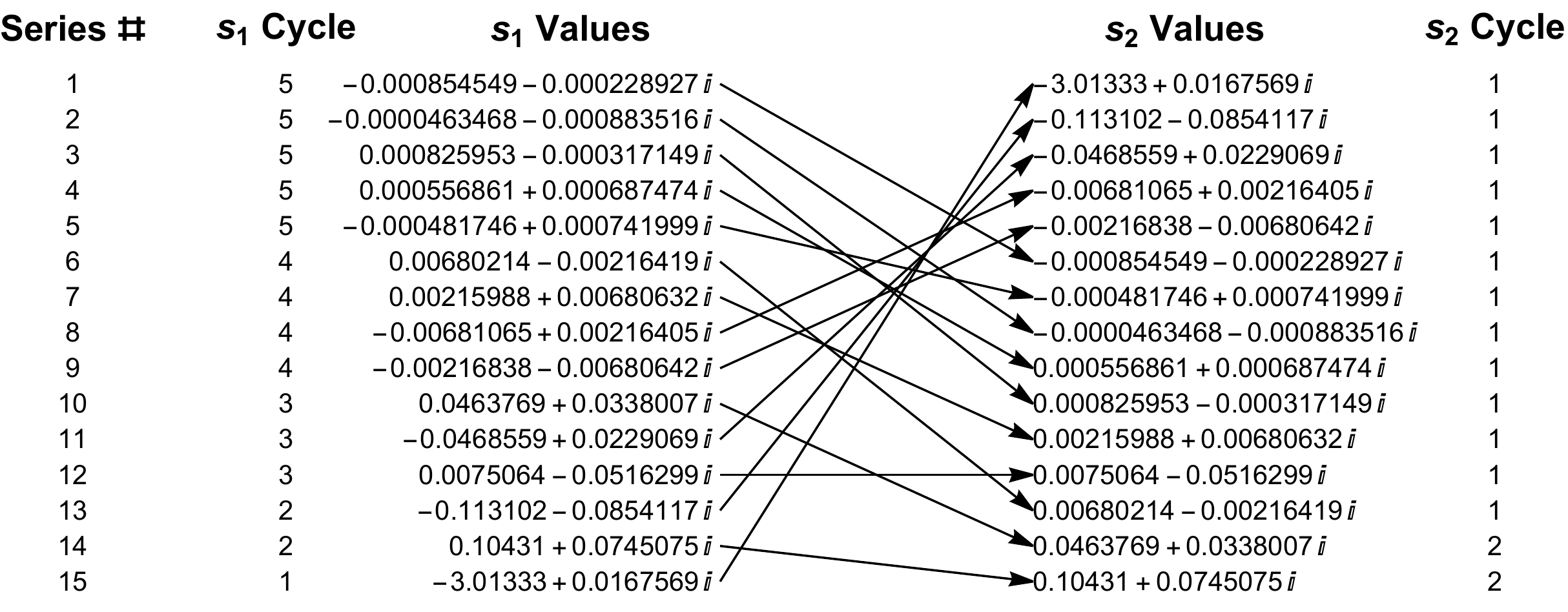}
\caption{Test Case $1$ analytic continuation between $s_1$ and $s_2$}
\label{figure:figure3}
\end{figure}

Figure \ref{figure:figure3} graphically illustrates the analytic continuation of all branch sheets at the expansion center $s_1$, into the branch sheets of $s_2$.  For example, the first series at $s_1$ is a $5$-cycle with the value at point $D$ of $-0.000855-0.000229i$.  Following the arrow from this point into the list of $s_2$ values, it continues onto sheet $6$ of $s_2$ which is a $1$-cycle.  The remaining four sheets of this $5$-cycle also continue onto $1$-cycles of $s_2$.  Therefore, the $5$-cycle at $s_1$ is analytically-continuous over $s_2$.  

 Consider now series $10$ of $s_1$ which is a $3$-cycle branch.  It continues onto series $14$ of $s_2$ which is a $2$-cycle branch.  Therefore, the $3$-cycle of $s_1$ is not analytically continuous over $s_2$ and so $s_2$ is the CLSP for this branch.  Likewise the $2$-cycle of series $14$ of $s_1$ continues onto the second sheet of the $2$-cycle of $s_2$ and therefore establishes the CLSP for this branch as well.  Therefore the CLSP for the $3$ and $2$-cycle branches is $s_2$. Next, analytic continuity of the $1$, $4$, and $5$-cycle branches are checked over $s_3$ and if any are continuous, further singular points are checked in this way until all CLSPs are found.   

%
%
%
%
\section{Creating accuracy profile functions}
\label{section:section005}

 Once the radii of convergence of a branch generator is determined, the accuracy of the series can be studied as a function of $|z|<R$ and series order $\mathscr{O}$. The accuracy of a series depends on four factors:

\begin{enumerate}
	\item \textbf{Precision of series:}
  	The precision of the series is limited to the precision of the expansion center, i.e., the singular points.  For the test cases below, singular points were computed to about $1000$ digits of precision.  Also, the precision of the associated power expansions decrease after each iteration of the Newton iteration step as shown in Figure \ref{figure:figure14}.  For example, terms of the $T$ series of Test Case 1 gradually drop to $918$ digits of precision near the end of $1000$ terms.  This trend limits the maximum accuracy of a particular value of the series to the precision of the terms used.  
	\item \textbf{Number of terms:}  The accuracy of an expansion was found to have a linear relation to the number of terms used in the expansion as illustrated in Figure \ref{fig112:subfig}.
	\item \textbf{The absolute value of $z$ relative to $R$:}  The accuracy of an expansion exhibited a logarithmic dependence on $|z|<R$ as shown in Figure \ref{fig110:subfig}.
	\item \textbf{Presence of nearby singular points}.  The accuracy of a series at $z$ is affected by the presence of nearby impinging singular points.  This is further explained below.
	\end{enumerate}

 An accuracy function $A(r_f,o)$ gives an expected accuracy of a series as a function of the radial ratio $0< r_f < 1$ and order, $o$, of a series.  An order function $O(r_f,e_a)$, for a given $r_f$ and desired accuracy $e_a$, returns the estimated order of a series needed for the desired accuracy.  The series is then be searched for the term $m$ corresponding to this order, and terms $1$ through $m$ used to evaluate the accuracy of the series.  Generator series are evaluated in their convergence domains and compared to more precise values of the function and fitted to accuracy functions.  Letting the expected accuracy $e_a=A(r_f,o)$, and solving for $o$ gives the order function.

\begin{figure}[ht] 
\centering
\includegraphics[scale=0.75]{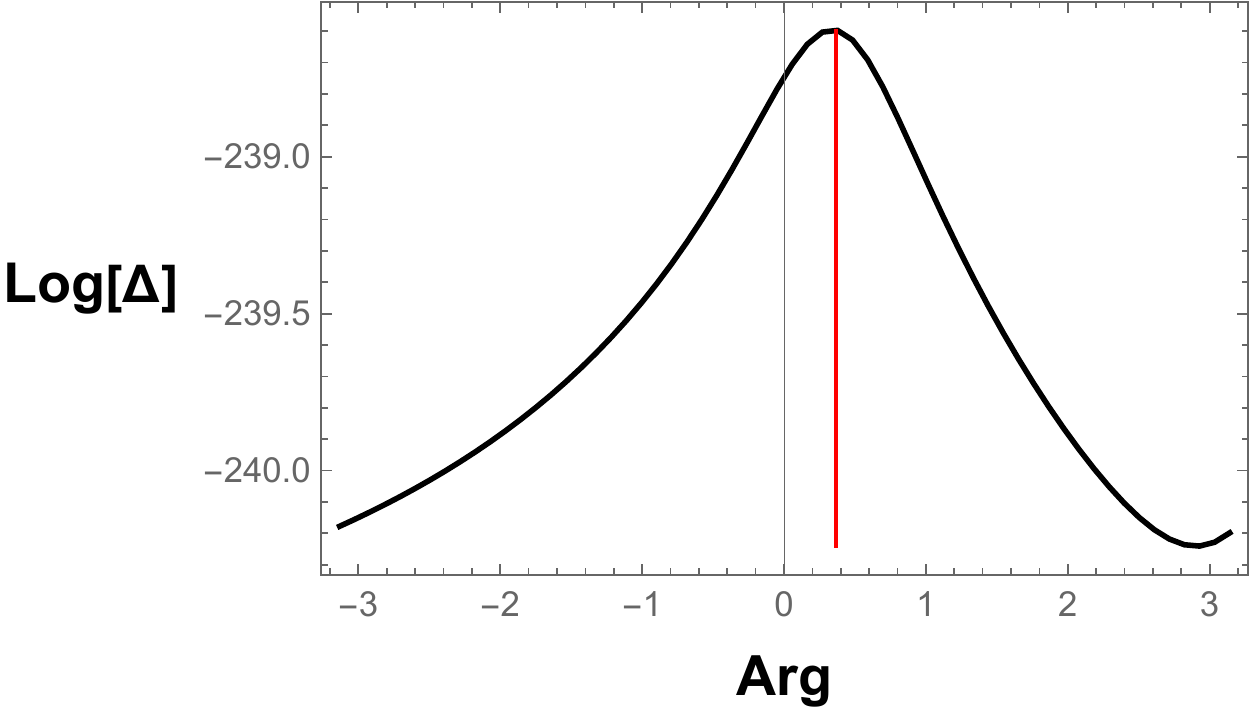}
\caption{Minimum arc accuracy diagram}
\label{figure:figure8}
\end{figure}

\begin{figure}[ht]
\centering
\subfloat[Accuracy points]{
\label{fig110:subfig:a} 
\includegraphics[scale=.7]{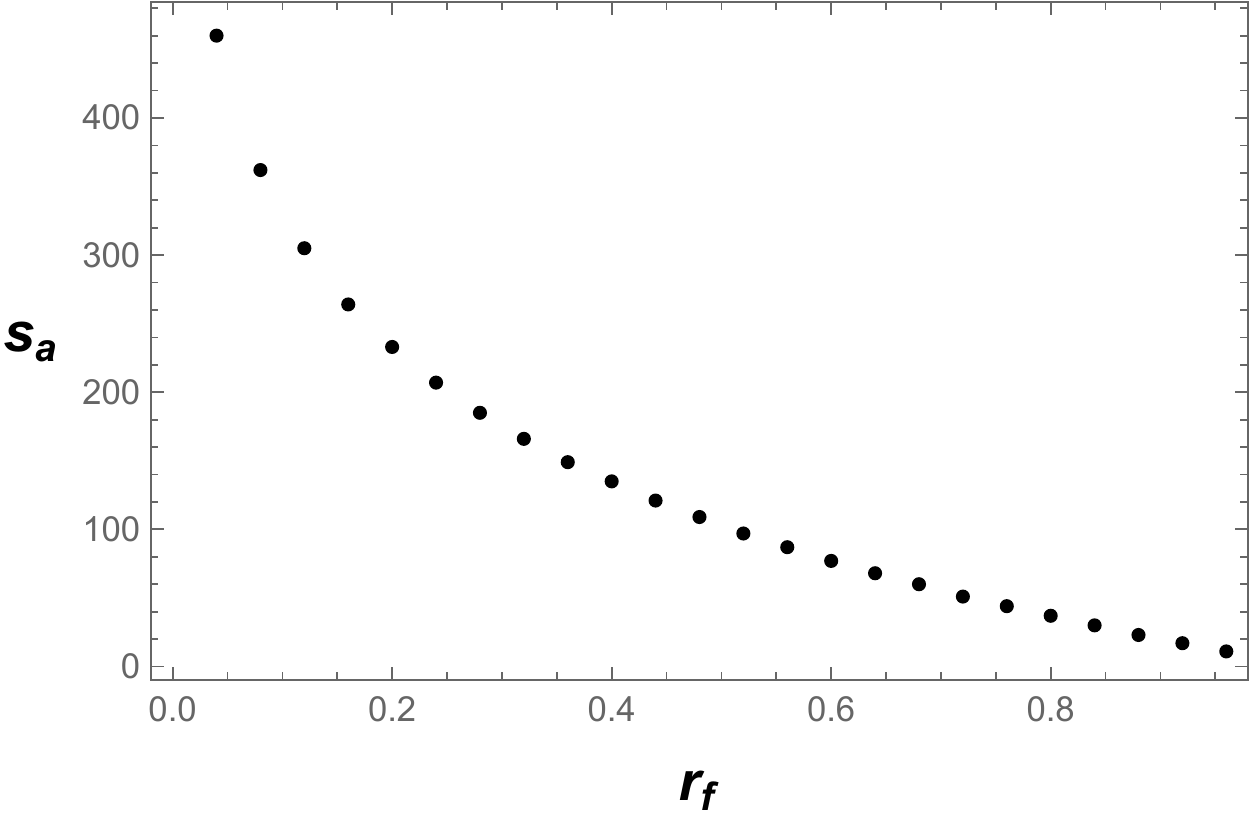}}
\hspace{.2in}
\subfloat[Fit function over accuracy points]{
\label{fig110:subfig:b} 
\includegraphics[scale=.7]{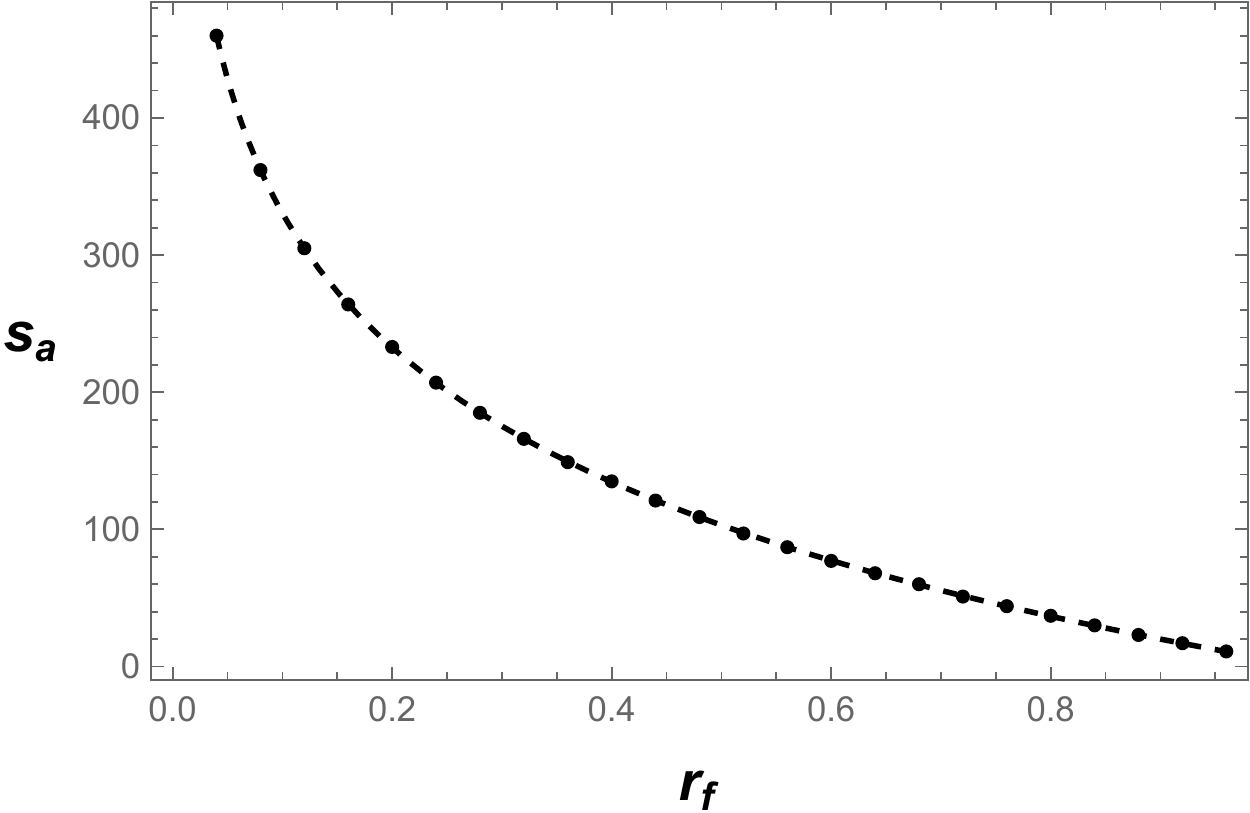}}
\hspace{.2in}
\caption{$F_5^{16}$ accuracy profile as a function of $r$, $1549$ terms}
\label{fig110:subfig}
\end{figure}

\begin{figure}[ht]
\centering
\subfloat[Accuracy points]{
\label{fig112:subfig:a} 
\includegraphics[scale=.7]{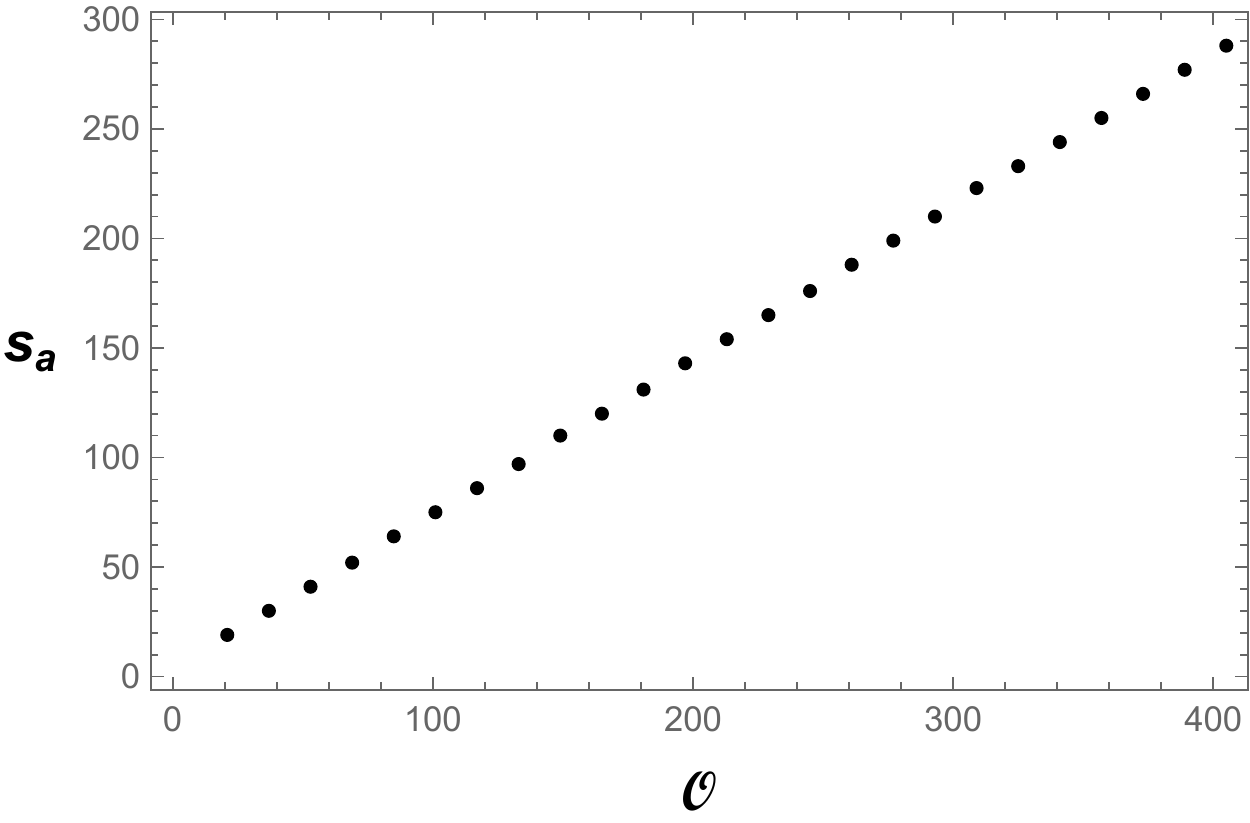}}
\hspace{.2in}
\subfloat[Fit function over accuracy points]{
\label{fig112:subfig:b} 
\includegraphics[scale=.7]{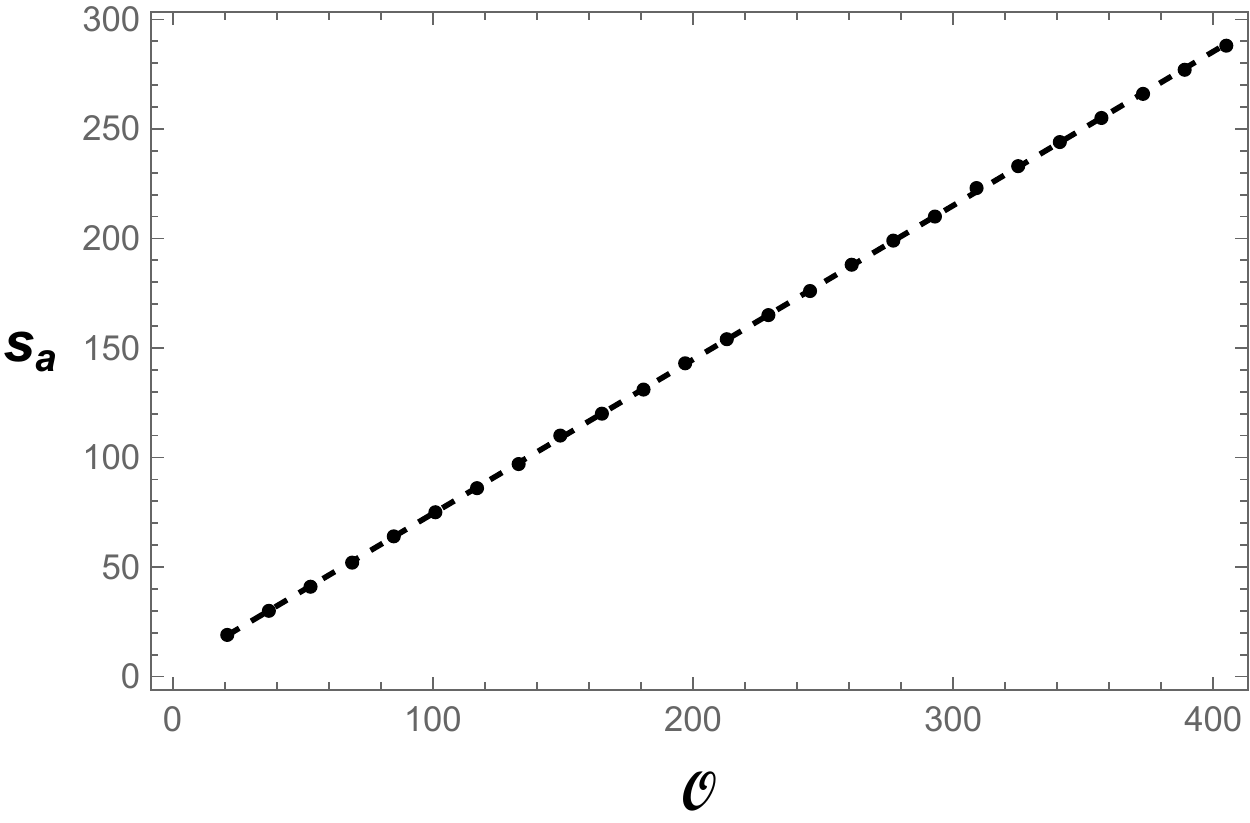}}
\hspace{.2in}
\caption{$F_5^{16}$ accuracy profile as a function of order with $r_f=1/5$}
\label{fig112:subfig}
\end{figure}

However, the accuracy of a power expansion is not constant along a circle $z=r_f R e^{it}$ but varies slightly along the outermost region of convergence depending on the presence of impinging singular points. The $F_5^{16}$ branch of Test Case $1$ has conjugate series $(1,2,3,4,5)$  with $R=s_{27}\approx 0.6413$.  The  impinging singular point of series $5$ is $s_{27}$.  Figure \ref{figure:figure8} is a plot of the log of the difference between the actual value of sheet $5$ and $1987$ terms of the series along a circle with $r_f=15/26$, that is $r=15/26 R$.  Notice how the difference is not constant around the circle but varies approaching a minimum accuracy at around an argument of $0.4$.  The red line in the plot is the argument of this branch sheet's ISP.  Notice the peak matches this line at $\log(\Delta)=-238.59$.  This and other cases suggest accuracy is affected by nearby impinging singular points.
The ISPs accuracy effect is however small.  Figure \ref{figure:figure8} reflects the log of the difference between the actual branch value and the series value which is more easily visualized.  The actual minimum and maximum difference is $(4.6\times 10^{-105},2.39\times 10^{-104})$.  This study does not include this variation in accuracy, rather, random points along the circle $z=r_f R e^{it}$ are used as test points although including this effect would improve the accuracy of the fit functions.

 Consider the accuracy results of $1287$ terms of this branch for $1/25\leq r_f \leq 24/25$.  Results of this test are shown in Figure \ref{fig110:subfig:a} and clearly shows a logarithmic trend.  
To this extent, the data points in Figure \ref{fig110:subfig:a} are fitted to $L(r)=a+b \log(r)$.  Figure \ref{fig110:subfig:b} shows this fit as the dashed black line $L(r)=3.3333 - 118.986\log(x)$.

Consider next the accuracy trend at $r_f=1/5$ for $20\leq o \leq 400$ shown in Figure \ref{fig112:subfig:a}.   The accuracy data follows a linear trend and is fitted to $G(o)=9.55371+0.0931429 o$ in Figure \ref{fig112:subfig:b}.

These accuracy trends are observed for all test cases studied in this paper.  
\begin{figure}[ht]
\centering
\subfloat[Accuracy points]{
\label{fig113:subfig:a} 
\includegraphics[scale=.6]{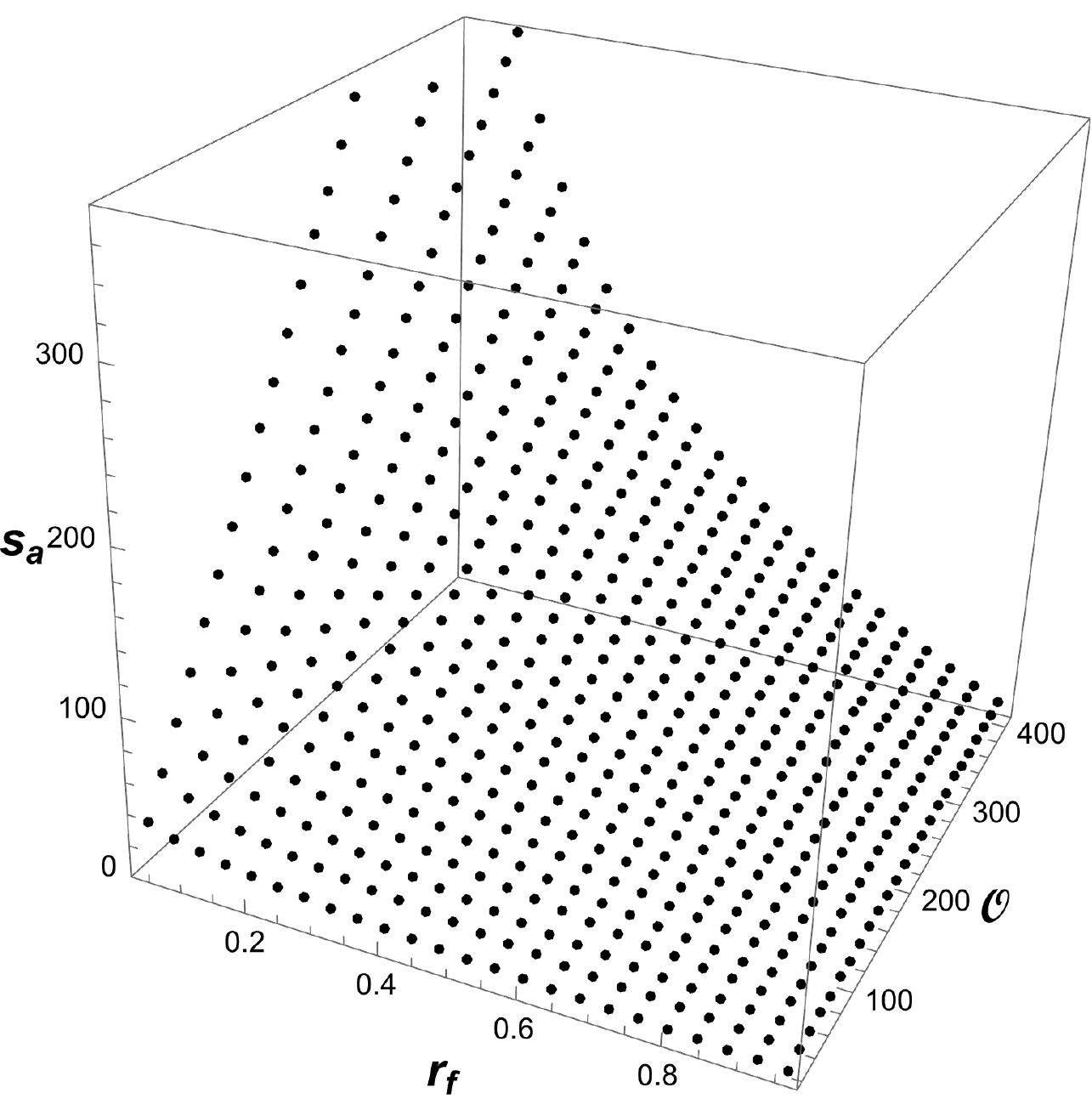}}
\hspace{.2in}
\subfloat[NonlinearModel fit over accuracy points]{
\label{fig113:subfig:b} 
\includegraphics[scale=.6]{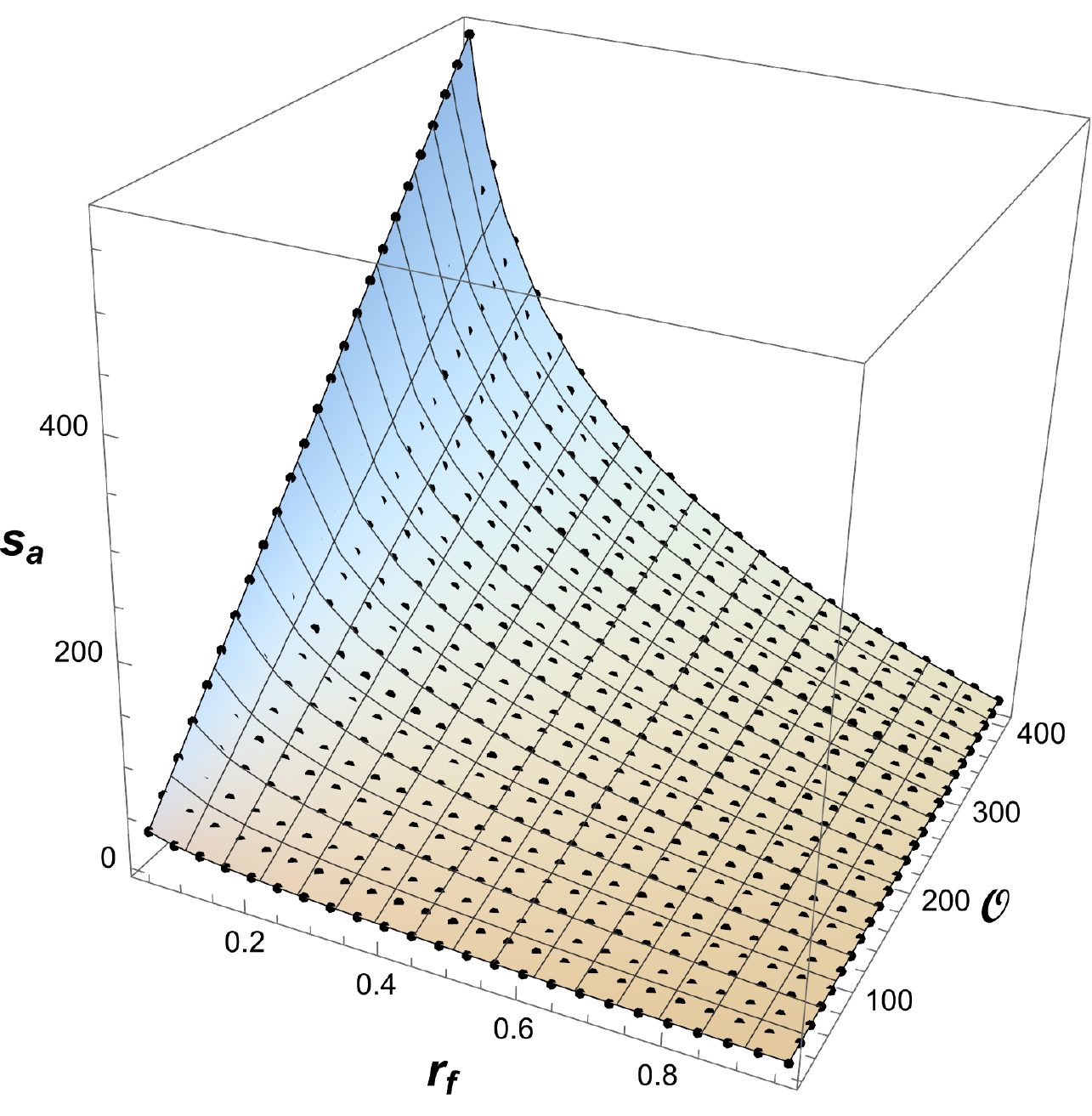}}
\hspace{.2in}
\caption{$F_5^{16}$ accuracy profile as a function of $r_f$ and $\mathscr{O}$}
\label{fig113:subfig}
\end{figure}

 When the accuracy data is plotted in space, the set of points shown in Figure \ref{fig113:subfig:a} are obtained.  Note in this figure the logarithmic trend over $r_f$ and the linear trend over the series order $\mathscr{O}$.  So it is reasonable to construct a log-linear accuracy function given by 
	\begin{equation}
	A(r_f,o)=a+b\log(r_f)+o(c+d \log(r_f))
	\label{eqn100}
	\end{equation}
	with $r_f$ the radial ratio and $o$, the order of the series using Mathematica's \texttt{NonlinearModelFit} function.  When this is done, the surface shown in \ref{fig113:subfig:b} with $A(r_f,o)=1.58647 -6.67332 \log(r)+o (0.000684195 -0.0871573 \log(r))$ is generated which also shows the accuracy data as black points superimposed on the surface.
	
 Letting $e_a=A(r_f,o)$ and solving for $o$ gives the order function
\begin{equation}
O(r_f,e_a)=\left\lceil\frac{e_a-a-b\log(r_f)}{c+d\log(r_f)}\right\rceil
\end{equation}
which estimates the order needed to achieve a particular accuracy $e_a$ at $r_f$.   
\section{Computing the Riemann surface genus using the Riemann-Hurwitz sum}
\label{sec:sec500}
 Once initial segments are computed for all singular points including the point at infinity, the Riemann surface genus is easily computed using the Riemann-Hurwitz formula in the form of 
\begin{equation}
\begin{aligned}
\mathscr{G}&=1+1/2\sum_{S}\sum_{i=1}^T(c_i-1)-D\\
&=1+1/2\; \mathscr{K}-D
\end{aligned}
\end{equation}

where $S$ is the set of singular points, and the double sum is over all conjugate groups of singular points with $c_i$ the cycle size.  $D$ is the degree of the function in $w$.  Note the Riemann-Hurwitz sum $\mathscr{K}$ must be an even number and serves as a cycle check-sum: if it's odd a cycle error has occurred.  However this is only a necessary condition for an error-free function ramification profile. As an example, Test Case $1$ with $179$ singular points has the ramification profile shown in Table \ref{table:table47b}.  A total of $174$ singular points each minimally ramify contributing a total of $174$ to $\mathscr{K}$.  $s_1$ contributes $10$, $s_{110}$ and $s_{111}$ contribute $16$ giving $\mathscr{K}=200$.  Then $\mathscr{G}=200/2+1-15=86$.  The time required for this calculation is predominantly the time to compute the singular points and initial segments.  
\section{Test results}

The following test cases include three data sets for each function:
\begin{enumerate}
	\item Timing data to compute the singular points, segments, power expansions, and CLSPs,
	\item Summary report of radii of convergence and accuracy results at a selected singular point,
	\item Ramification profile summarizing the cycle geometry at all singular points.
\end{enumerate}
%
%
%
%
\subsection{Test Case 1: $15$-degree function expanded at the origin,$\mathscr{G}=86$}

 \label{testCase1}

\begin{equation}
f_1(z,w)=\left(z^{30}+z^{32}\right)+\left(z^{14}+z^{20}\right) w^5+\left(z^5+z^9\right) w^9+\left(z+z^3\right) w^{12}+6 w^{14}+\left(2+z^2\right) w^{15}
	\label{eqn:eqn409}
	\end{equation}

 Timing data is shown in Table \ref{table:table44b}.  The base expansions are done to at least $1000$ terms and the comparison series to at least $100$ terms at $1000$ digits of precision.  
\begin{table}[ht]
\caption{Test Case $1$ Timing Data}
$
\begin{array}{c}
\begin{array}{|c|c|c|c|c|c|}
\hline
 \begin{array}{c}
 \text{Singular} \\
  \text{points}
	\end{array}&
	\begin{array}{c}
 \text{Initial} \\
  \text{segments}
	\end{array}&
	\begin{array}{c}
 \text{Base gen.} \\
  \text{expansions}
	\end{array}&
	\begin{array}{c}
 \text{Comparison} \\
  \text{expansions}
	\end{array}&
	\text{CT} &
	\text{IT}
	\\
\hline
 $(179,1.9\,s)$ & $8\,s$ &$(5,13.4\,m)$ & $(125,37.5\,m)$ & $1.3\,m$ & $2.7\,m$\\
 \hline 
\end{array} \\
\scriptsize \textbf{Singular points: }\text{(Total points, time), }\textbf{Initial segments: }\text{time }\\
\scriptsize\textbf{Base expansion: }\text{(Total generators,time), }\textbf{Comparison expansions: }\text{(Total sing. pts.,time)}\\
		\scriptsize\textbf{CT: }\text{comparison test time, }\textbf{IT: }\text{Integration Test time}\\
		\scriptsize\textbf{s: }\text{Seconds, }\textbf{m: }\text{Minutes }\textbf{h: }\text{hours}
\end{array}
$
\label{table:table44b}
\end{table}

Table \ref{table:table47} summarizes the CLSP and accuracy results. 

\begin{table}[ht]
\caption{Test Case $1$ Summary Report at the origin}
$
\begin{array}{c}
\begin{array}{|c|c|c|c|c|c|c|c|c|}
\hline
 \text{Type} & \text{CLSP} & \text{R} & \text{Terms} & \text{a} & \text{b} & \text{c} & \text{d} & \text{Var} \\
\hline
 \, F_5^{16} & 27 & 0.641328 & 1987 & 3.19696 & -0.329926 & 0.00567716 & -0.433846 & 0.23481 \\
 \, F_4^9 & 7 & 0.504901 & 2022 & 3.4831 & -0.219358 & 0.00355264 & -0.434293 & 0.228403 \\
 \, F_3^4 & 2 & 0.166817 & 1017 & 3.42982 & -0.287641 & 0.00466601 & -0.434467 & 0.173804 \\
 \, V_2^1 & 2 & 0.166817 & 1024 & 3.5265 & -0.224677 & 0.00418632 & -0.434702 & 0.259237 \\
 \text{T} & 118 & 1.09352 & 1024 & 2.83953 & -0.287592 & 0.00193844 & -0.434512 & 0.565284 \\
\hline
\end{array}
\\
 \scriptsize \textbf{Type: }\text{Branch type (see Appendix \ref{appendix:appA}), }\small \textbf{CLSP: }\text{Convergence-limiting singular point, }\\
\scriptsize\textbf{R: }\text{Radius of convergence, }\small \textbf{Terms: }\text{Series length, }
    \scriptsize \textbf{(a,b,c,d): }\text{constants for $A(r_f,o)$, } \textbf{Var: }\text{Variance}\\
		\end{array}
	$
\label{table:table47}
\end{table}

The accuracy constants $(a,b,c,d)$ can be used to determine the approximate order needed to achieve a desired accuracy.  For example, if a value of the $F_5^{16}$ branch at $z=1/3 R e^{3 \pi i/4}$ to $20$ digits of accuracy is desired, that is, $e_a=20$, solve $O(1/3,20)$ where
$$
O(r_f,e_a)=\left\lceil \frac{a+0.329926 \log (r)-3.19696}{0.00567716\, -0.433846 \log (r)}\right\rceil.
$$
This gives an order estimate of $35$.  A simple scan of the expansion exponents can determine the series term corresponding to an order of $35$.  This turns out to be term $99$.  The first $99$ terms of the generator series at $z$ leads to a value accurate to approximately $20$ digits.  In this case, the accuracy is $20$ digits as shown by comparing the value to the corresponding root of $f_1(z,w)=0$.  The estimated accuracy however will often differ from the actual error by a small amount as shown by the variance in Table \ref{table:table47}.

 The constants $(a,b,c,d)$ for each branch in Table \ref{table:table47} are similar in size.  This means that a single average accuracy function could be generated for all the branches in this case.  However in other cases, the accuracy constants can be quite different.

 Table \ref{table:table47b} summarizes the ramification profile describing the cycle geometry at all singular points.  The set of singular points minimally-ramified is $\{\bar{s}_{n}\}$ with ramification $(2,[13,1])$ signifying a single $2$-cycle branch and $13$ single-cycle branches.  Singular points with higher ramifications are listed separately or $\{u_i\}_n$ for multiple singular points with the same ramification.

\begin{table}[ht]
\caption{Test Case $1$ Ramification profile, $\mathscr{K}=200$}
$
\small
\begin{array}{c}
\begin{array}{|c|c|}
\hline
 \text{Singular point} & \text{Cycles} \\
\hline
 s_1 & (1,2,3,4,5) \\
s_{110} & (9,[6,1]) \\
s_{111} & (9,[6,1]) \\
\{p_i\}_2 & [15,1]\\
\{\bar{s}\}_{174} &      (2,[13,1])\\
s_{\infty} & \text{[15,1]}\\
\hline
\end{array}
\\
\scriptsize \textbf{$\{p_i\}_2: $}\text{ set of two poles }\\
 \scriptsize \textbf{[m,n]: }\text{$n$-cycle branches, $m$ total}\\
  	\scriptsize \textbf{$\{\bar{s}\}_n: $}\text{ remaining $n$ singular points minimally ramified}\\
		\end{array}
	$
\label{table:table47b}
	
\end{table}

%
%
%
%
\subsection{Test Case $2$: $4$-degree function with polynomial solution at the origin, $\mathscr{G}=0$}\label{testCase2}
\begin{equation}
f_2(z,w)=\left(1-3 z+3 z^2-1 z^3\right)+\left(-4+8 z-4 z^2\right) w+(6-6 z) w^2+(-4) w^3+(1) w^4
		\label{eqn:eqn003b}
\end{equation}

This function has two singular points $\{0,1\}$ and a $4$-cycle polynomial solution at the origin:
\begin{equation}
\begin{aligned}
w_1(z)&=1-z^{1/4}+z^{1/2}-z^{3/4}\\
w_2(z)&=1-iz^{1/4}-z^{1/2}+iz^{3/4}\\
w_3(z)&=1+z^{1/4}+z^{1/2}+z^{3/4} \\
w_4(z)&=1+iz^{1/4}-z^{1/2}-iz^{3/4}.
\end{aligned}
\end{equation}

 In this case, the modular operation of the Newton iteration step (\ref{eqn:eqn989}) returns a finite polynomial after reaching the maximum zero modular value $N_{zm}$.  And since the solutions are finite, $R=\infty$.  However this can only be true if the function does not ramify or is non-polar at $z=1$.  That is, the singularity at $z=1$ is removable.  We can show this as follows:  

The expansions at $z=1$ are four $1$-cycles with starting terms:
\begin{equation}
\begin{aligned}
w_1(z)&=-0.5 z+0.0625 z^2-0.03125 z^3+0.0205078 z^4+\cdots\\
w_2(z)&=(-0.5-0.5 I) z+(0.125 +0. I) z^2-(0.0625 -0.015625 I) z^3+\cdots\\
w_3(z)&=(-0.5+0.5 I) z+(0.125 +0. I) z^2-(0.0625 +0.015625 I) z^3+\cdots \\
w_4(z)&=4+1.5 z-0.3125 z^2+0.15625 z^3+\cdots.
\label{eqn:eqn602b}
\end{aligned}
\end{equation}
The solutions to $f_2(1,w)=0$ are $\{0,0,0,4\}$ and these are the values of (\ref{eqn:eqn602b}) at $z=0$ and note 
\begin{equation}
\begin{aligned}
\frac{dw}{dz}&=\lim_{(z,w)\to(1,0)} \left(-\frac{f_z(z,w)}{f_w(z,w)}\right)\to\frac{0}{0}\\
\frac{dw}{dz}&=\lim_{(z,w)\to(1,4)} \left(-\frac{f_z(z,w)}{f_w(z,w)}\right)=-\frac{3}{2}.
\end{aligned}
\end{equation}

However the derivatives of (\ref{eqn:eqn602b}) at $z_r=0$ are finite so that we can immediately solve for the indeterminate limits
\begin{equation}
\frac{dw}{dz}=\lim_{(z,w)\to(1,0)} \left(-\frac{f_z(z,w)}{f_w(z,w)}\right)=\{-0.5,-0.5-0.5i,-0.5+0.5i\}.
\end{equation}

And since the function fully-ramifies at the origin, the four expansions at $z_a=1$ all have $R=1$.  This is confirmed by the Root Test and both series comparison and integration tests.  Accuracy results are given in Table \ref{table:table52}.  And by virtue of the singularity at the origin, the expansions at $z=1$ all have radii of convergence of $1$.  Finally, as this branch is equivalent to the function $f(z)=z^{1/4}$, the ramification at infinity will also be $4$-cycle and thus the genus is $(3+3)/2+1-4=0$.  Timing for this case was minimal.
\begin{table}[ht]
\caption{Test Case $2$ Summary Report at $s_2=1$}
$
\small
\begin{array}{|c|c|c|c|c|c|c|c|c|}
\hline
 \text{Type} & \text{CLSP} & \text{R} & \text{Terms} & \text{a} & \text{b} & \text{c} & \text{d} & \text{Var} \\
\hline
 \text{1E} & 1 & 1. & 1024 & 3.07174 & -0.302067 & 0.00142759 & -0.434367 & 0.156447 \\
 \text{2E} & 1 & 1. & 1024 & 2.83926 & -0.321447 & 0.00161223 & -0.434362 & 0.164684 \\
 \text{3E} & 1 & 1. & 1024 & 2.83926 & -0.321447 & 0.00161223 & -0.434362 & 0.164684 \\
 \text{T} & 1 & 1. & 1024 & 2.59935 & -0.324381 & 0.00175722 & -0.434358 & 0.195594 \\
\hline
\end{array}
	$
\label{table:table52}
\end{table}

%
%
%
%
\subsection{Test Case $3$:  $34$-degree function expanded at $s_{487}$, $\mathscr{G}=264$}\label{testCase3}
 
\begin{equation}
\begin{aligned}
f_3(z,w)&=\left(-\frac{1043}{60}-\frac{5}{3} z^2+2 z^3-4 z^4-\frac{6}{5} z^5-\frac{2}{3} z^9+\frac{8}{3} z^{14}+\frac{25}{4} z^{15}+4 z^{16}\right)\\
&+\left(\frac{11}{3}\right) w^5+\left(-\frac{8}{3}\right) w^{12}+\left(-\frac{38}{5}-\frac{1}{2} z\right)
   w^{34}\\
	&=0
	\end{aligned}
	\label{eqn:eqn1095}
	\end{equation}
This function was selected to stress-test the methods of finding CLSPs and also to study the fully-ramified branch at infinity.  The function has $493$ singular points with $s_{487}$ ramifying into a $22$-cycle and $12$ single cycles and having $13$ nearest neighbors with an average separation distance of $10^{-35}$.

Root Tests results were not precise enough to distinguish individual singular points near the expansion center and returned $1.056\times 10^{-35}$ as the estimate for the single-cycle branches and $1.09\times 10^{-35}$ for the $22$-cycle branch identifying $s_{487}$ as the estimated CLSP of the base expansions.  This in itself is not an error but rather attempting to use the Root Test at an extremely small tolerance.  The comparison sequence however can be set to a reasonable number of singular points and in this case was set to the first $25$ singular point in the singular sequence.  

However, the comparison test failed to identify CLSPs at the default settings as the base series were not accurate enough to meet the comparison tolerances of $10^{-6}$.  The Integration Test likewise failed to identify CLSPs with a working precision of $40$ and default integration method.    However, increasing the integration working precision of the integration test from $40$ to $80$ with ``StiffnessSwitching'' method met the tolerances and found all CLSPs in $11.5$ minutes.   These results are shown in Table \ref{table:table29b}.  
\begin{table}[h]
\caption{Test Case $3$ Timing Data}
$
\begin{array}{|c|c|c|c|c|c|}
\hline
 \begin{array}{c}
 \text{Singular} \\
  \text{points}
	\end{array}&
	\begin{array}{c}
 \text{Initial} \\
  \text{segments}
	\end{array}&
	\begin{array}{c}
 \text{Base} \\
  \text{expansion}
	\end{array}&
	\begin{array}{c}
 \text{Comparison} \\
  \text{expansions}
	\end{array}&
	\text{CT} &
	\text{IT}
	\\
\hline
 $(493,54\,s)$ & $(50\,s)$ &$(13,\,13.2m)$ & $(25,\,7.5m\.h)$ & $-$ & $11.5 m$  \\
 \hline 
\end{array} 
$
\label{table:table47d}
\end{table} 
\begin{table} [h]
\caption{Test Case $3$ Summary Report at $s_{487}$}
$
\small
\begin{array}{c}
\begin{array}{|c|c|c|c|c|c|c|c|c|}
\hline
 \text{Type} & \text{CLSP} & \text{R} & \text{Terms} & \text{a} & \text{b} & \text{c} & \text{d} & \text{Var} \\
\hline
 \text{1T} & 493 & 1.0511\times 10^{-35} & 1024 & 2.57488 & -0.408571 & 0.00198314 & -0.434325 & 0.246471 \\
 \text{2T} & 488 & 1.0511\times 10^{-35} & 1024 & 2.43652 & -0.456495 & 0.00189591 & -0.434416 & 0.221283 \\
\text{3T} & 489 & 1.0511\times 10^{-35} & 1024 & 2.42991 & -0.493089 & 0.00181696 & -0.434467 & 0.24568 \\
 \text{4T} & 485 &1.0511\times 10^{-35} & 1024 & 2.58246 & -0.391823 & 0.00188475 & -0.434426 & 0.304147 \\
 \text{5T} & 486 & 1.0511\times 10^{-35} & 1024 & 2.49944 & -0.447571 & 0.00205355 & -0.434356 & 0.244369 \\
 \text{6T} & 482 & 1.0511\times 10^{-35} & 1024 & 2.60673 & -0.384546 & 0.0019365 & -0.434366 & 0.241232 \\
 \text{7T} & 481 & 1.0511\times 10^{-35} & 1024 & 2.52156 & -0.458827 & 0.00198745 & -0.434316 & 0.238367 \\
 \text{8T} & 484 & 1.0511\times 10^{-35} & 1024 & 2.64395 & -0.365494 & 0.00186476 & -0.434422 & 0.242922 \\
 \text{9T} & 483 & 1.0511\times 10^{-35} & 1024 & 2.61553 & -0.379034 & 0.00196351 & -0.434371 & 0.247183 \\
 \text{10T} & 491 & 1.0511\times 10^{-35} & 1024 & 2.66054 & -0.378256 & 0.00185278 & -0.434409 & 0.229304 \\
 \text{11T} & 490 & 1.0511\times 10^{-35} & 1024 & 2.56919 & -0.427104 & 0.00188795 & -0.434372 & 0.248727 \\
 \text{12T} & 492 & 1.0511\times 10^{-35} & 1024 & 2.60019 & -0.386299 & 0.00190533 & -0.434368 & 0.246803 \\
 \, P_{22}^{-1} & 492 & 1.0511\times 10^{-35} & 2021 & 1.64383 & -0.210138 & 0.0150419 & -0.433515 & 0.183742 \\
\hline
\end{array}\\
\\
\small \textbf{Actual values of $R$ differ from the $1.0511\times 10^{-35}$ value given for brevity}\\
\small \textbf{For precise values of $R$ use $R=|s_{487}-s_{CLSP}|$} \\
\end{array}
$
\label{table:table29b}
\end{table}

 The function has a fully-ramified $34$-cycle branch at infinity which means the CLSP is $s_2$ relative to infinity.  The default setting for the comparison test were not sufficient to identify this CLSP. The integration test using default integration methods also failed but succeeded in identifying $s_2$ as the CLSP using Mathematica's \texttt{NDSolve} \texttt{StiffnessSwitching} method. 
%
%
\begin{table}[ht]
\caption{Test Case $3$ Ramification profile, $\mathscr{K}=594$}
$
\small
\begin{array}{c}
\begin{array}{|c|c|}
\hline
 \text{Singular point} & \text{Cycles} \\
\hline
 \{u_i\}_{16} & (5,[29,1])) \\
s_{487} & (22,[12,1]) \\
\{\bar{s}\}_{476} &      (2,[32,1])\\
s_{\infty} & (34)\\
\hline
\end{array}
\\
\tiny \textbf{$\{u_i\}_n$}\text{=$\{51,52,57,58,134,177,178,255,256,281,282,357,358,426,427,440\}$}\\
		\end{array}
	$
\label{table:table49c}
	
\end{table}

%
%
%
%
\subsection{Test Case 4:  $35$-degree function expanded at infinity,$\mathscr{G}=32$}
\label{testCase4}

\begin{equation}
\begin{aligned}
f_4(z,w)&=\left(-\frac{31}{10}+\frac{179}{30} z+\frac{1}{4} z^2\right)+\left(-\frac{7}{4}\right) w^2+(4)
   w^3+\left(-\frac{1}{2}-\frac{5}{2} z\right) w^8+\left(\frac{11}{3}\right) w^{10}+\left(6+\frac{5}{2} z\right) w^{14}\\
	&+(5)
   w^{18}+\left(-\frac{64}{15}\right) w^{20}+\left(\frac{11}{2}-\frac{1}{2} z^2\right) w^{22}+\left(-\frac{9}{2}+\frac{7}{3}
   z\right) w^{25}+\left(\frac{18}{5}-\frac{3}{4} z^2\right) w^{28}\\
	&+\left(-\frac{3}{2}-1 z\right)
   w^{33}+\left(-\frac{8}{3}\right) w^{35}\\
	&=0
	\end{aligned}
	\end{equation}

 First consider $f_4(z,w)$ which has $127$ finite singular points all of which are minimally-ramified.  In order to obtain the ramification at infinity, $f_4(z,w)$ is transformed to $\displaystyle g_4(z,w)=z^\delta f_4\left(\frac{1}{z},w\right)$ with $\delta$ the largest power of $z$ in $f_4(z,w)$.  In this case $\delta=2$ giving:
	\begin{equation}
	\begin{aligned}
	g_4(z,w)&=\left(\frac{1}{4}+\frac{179}{30} z-\frac{31}{10} z^2\right)+\left(-\frac{7}{4} z^2\right) w^2+\left(4 z^2\right)
   w^3+\left(-\frac{5}{2} z-\frac{1}{2} z^2\right) w^8+\left(\frac{11}{3} z^2\right) w^{10}\\
	&=+\left(\frac{5}{2} z+6 z^2\right)w^{14}+\left(5 z^2\right) w^{18}+\left(-\frac{64}{15} z^2\right) w^{20}+\left(-\frac{1}{2}+\frac{11}{2} z^2\right)w^{22}\\
	&+\left(\frac{7}{3} z-\frac{9}{2} z^2\right) w^{25}+\left(-\frac{3}{4}+\frac{18}{5} z^2\right)
   w^{28}+\left(-z-\frac{3}{2} z^2\right) w^{33}+\left(-\frac{8}{3} z^2\right) w^{35}\\
	&=0.
	\end{aligned}
	\end{equation}
	
	 Then an expansion of $w(z)$ at $z_a=0$ through $g_4(z,w)$ is the expansion of $w(z)$ defined by $f_4(z,w)$ at infinity.  The base series and comparison series were next computed relative to $g_4(z,w)$.  Table \ref{table:table147c} is the timing summary.
\begin{table}[h]
\caption{Test Case $4$ Timing Data at infinity}
$
\begin{array}{|c|c|c|c|c|c|}
\hline
 \begin{array}{c}
 \text{Singular} \\
  \text{points}
	\end{array}&
	\begin{array}{c}
 \text{Initial} \\
  \text{segments}
	\end{array}&
	\begin{array}{c}
 \text{Base} \\
  \text{expansion}
	\end{array}&
	\begin{array}{c}
 \text{Comparison} \\
  \text{expansions}
	\end{array}&
	\text{CT} &
	\text{IT}
	\\
\hline
 $(128,2\,s)$ & $(15\,s)$ &$(30,3.4\,h)$ & $(32,68\,m)$ & $2.6\,m$ & $7.4\,m$ \\
 \hline 
\end{array} 
$
\label{table:table147c}
\end{table}	

The expansion at infinity ramified as $(5,2,[28,1])$.  Table \ref{table:table294} is a Summary Report.  Consider the $5$-cycle branch expansions of $P_5^{-1}$ of $g_4(z,w)$ with $R=0.04874$.  These are expansions of $w(z)$ centered at infinity which means the five values of $w(z)$ associated with this branch at $\displaystyle z_a=\frac{1}{z_r}$ are the same five branch values $\{v_s\}$ at $z_r$.  For example, let $z_a=200+100i$ which is outside the domain of finite singular points of $f_4$ and therefore $z_r=\displaystyle\frac{1}{200+100i}=0.004-0.002i$ is closer to infinity than the nearest singular point of $f_4$.  Then the $P_5^{-1}$ expansions will converge at $z_r$.  The $35$ values of $w(z_a)$ are easily found by solving for the roots $\{w_i\}=f_4(z_a,w)=0$.  Among these roots are the five values of $P_5^{-1}$ at  $z_r$:
$$
\mscriptsize{-2.83979-0.469744 i,-1.3576+2.63902 i,-0.501188-2.89321 i,2.05904\, +2.02105 i,2.62398\, -1.28185 i.} 
$$

 The following are the roots $\{w_i\}$ with the above branch values highlighted in red. 
$$
\begin{aligned}
 &\mscriptsize{-3.92407+9.53357 i,\textcolor{red}{-2.83979-0.469744 i},\textcolor{red}{-1.3576+2.63902 i},-0.940335+0.00182458 i,-0.922871+0.235934 i,} \\
&\mscriptsize{-0.922161-0.232486 i,-0.888861-0.448904 i,-0.886022+0.452484 i,-0.783866-0.586343 i,-0.780327+0.586443 i,}\\
&\mscriptsize{-0.604223-0.72708 i,-0.601147+0.730563 i,\textcolor{red}{-0.501188-2.89321 i},-0.403707-0.850957 i,-0.399549+0.853659 i,}\\
&\mscriptsize{-0.187928-0.945784 i,-0.181958+0.946171 i,-0.000537808-0.999928 i,-0.000128033+1.00055 i,0.182496\, -0.945415 i,}\\
&\mscriptsize{0.187148\, +0.94525 i,0.399708\, -0.85302 i,0.403067\, +0.850813 i,0.601121\, -0.730076 i,0.603749\, +0.727117 i,}\\
&\mscriptsize{0.781083\, -0.585674 i,0.783141\, +0.585525 i,0.887975\, -0.453228 i,0.889622\, +0.446953 i,0.92263\, +0.231865 i,}\\
&\mscriptsize{0.923506\, -0.236407 i,0.940601\, -0.00209421 i,\textcolor{red}{2.05904\, +2.02105 i},\textcolor{red}{2.62398\, -1.28185 i},3.9374\, -9.5466 i}
\end{aligned}
$$
Another way to visualize the expansions at infinity is to consider the ramification of $f_4(z,w)$ outside a disc containing all the finite singular points, that is a disc $r<|s_{127}|\approx 27.15$.  This ramification is the same $(5,2,[28,1])$ ramification as that at infinity and is in fact the same set of branches.  That is, we can compute the ramification at infinity for $f_4(z,w)$ by simply computing the ramification of $f_4(z,w)$ at say $r=28$.  
See \href{https://arxiv.org/pdf/1901.03996.pdf}{On the branching geometry of algebraic functions} for a method of doing this.
\begin{table} [H]
\caption{Test Case $4$ summary report at the point of infinity, i.e., $g_4$ expanded at $z=0$}
$
\small
\begin{array}{c}
\begin{array}{|c|c|c|c|c|c|c|c|c|}
\hline
 \text{Type} & \text{CLSP} & \text{R} & \text{Terms} & \text{a} & \text{b} & \text{c} & \text{d} & \text{Var} \\
\hline
 \text{1T} & 23 & 0.0416879 & 1024 & 3.28983 & -0.360237 & 0.00187946 & -0.434348 & 0.186178 \\
 \text{2T} & 2 & 0.0368337 & 1024 & 3.45719 & -0.343809 & 0.00191956 & -0.434455 & 0.25136 \\
 \text{3T} & 3 & 0.0368337 & 1024 & 3.40354 & -0.405533 & 0.0019289 & -0.434433 & 0.257028 \\
 \text{4T} & 2 & 0.0368337 & 1024 & 3.34071 & -0.453646 & 0.00192991 & -0.434402 & 0.289739 \\
 \text{5T} & 3 & 0.0368337 & 1024 & 3.27402 & -0.498051 & 0.00198714 & -0.434375 & 0.335869 \\
 \text{6T} & 5 & 0.0369532 & 1024 & 3.67988 & -0.426256 & 0.00192521 & -0.43436 & 0.253784 \\
 \text{7T} & 4 & 0.0369532 & 1024 & 3.70872 & -0.420795 & 0.00196974 & -0.434302 & 0.250291 \\
 \text{8T} & 24 & 0.0419264 & 1024 & 3.26493 & -0.392033 & 0.001903 & -0.434354 & 0.220053 \\
 \text{9T} & 25 & 0.0419264 & 1024 & 3.44163 & -0.341589 & 0.00175748 & -0.434415 & 0.21876 \\
 \text{10T} & 13 & 0.0401948 & 1024 & 3.22901 & -0.429983 & 0.00206575 & -0.434317 & 0.270064 \\
 \text{11T} & 12 & 0.0401948 & 1024 & 3.16109 & -0.464687 & 0.00199673 & -0.43438 & 0.303588 \\
 \text{12T} & 13 & 0.0401948 & 1024 & 3.17537 & -0.444561 & 0.00193504 & -0.434374 & 0.239211 \\
 \text{13T} & 12 & 0.0401948 & 1024 & 3.34242 & -0.380652 & 0.00192249 & -0.43438 & 0.221718 \\
 \text{14T} & 28 & 0.0590922 & 1024 & 3.85445 & -0.402529 & 0.00188963 & -0.434352 & 0.247296 \\
 \text{15T} & 29 & 0.0590922 & 1024 & 3.82279 & -0.484877 & 0.00188897 & -0.434324 & 0.247753 \\
 \text{16T} & 10 & 0.0391605 & 1024 & 3.3216 & -0.409763 & 0.00199395 & -0.434323 & 0.220968 \\
 \text{17T} & 11 & 0.0391605 & 1024 & 3.32488 & -0.391127 & 0.00195106 & -0.434362 & 0.228586 \\
 \text{18T} & 10 & 0.0391605 & 1024 & 3.24861 & -0.444286 & 0.00197163 & -0.434382 & 0.270096 \\
 \text{19T} & 11 & 0.0391605 & 1024 & 3.265 & -0.437944 & 0.00201653 & -0.434335 & 0.270439 \\
 \text{20T} & 16 & 0.0408681 & 1024 & 3.33897 & -0.353798 & 0.00190399 & -0.434374 & 0.202345 \\
 \text{21T} & 15 & 0.0408681 & 1024 & 3.28935 & -0.354593 & 0.00189579 & -0.43441 & 0.207717 \\
 \text{22T} & 8 & 0.0383842 & 1024 & 3.83458 & -0.418136 & 0.00191619 & -0.434332 & 0.237896 \\
 \text{23T} & 9 & 0.0383842 & 1024 & 3.93489 & -0.370628 & 0.00194949 & -0.434367 & 0.218989 \\
 \text{24T} & 7 & 0.0378856 & 1024 & 3.36544 & -0.388662 & 0.0019894 & -0.434436 & 0.264613 \\
 \text{25T} & 6 & 0.0378856 & 1024 & 3.35258 & -0.425374 & 0.00204086 & -0.434361 & 0.309241 \\
 \text{26T} & 7 & 0.0378856 & 1024 & 3.3512 & -0.423055 & 0.00192371 & -0.434414 & 0.251619 \\
 \text{27T} & 6 & 0.0378856 & 1024 & 3.43569 & -0.376695 & 0.00187883 & -0.434427 & 0.252273 \\
 \text{28T} & 26 & 0.0427618 & 1024 & 3.29122 & -0.404697 & 0.00182958 & -0.434401 & 0.215776 \\
 \, P_5^{-1} & 27 & 0.0487394 & 1020 & 3.27798 & -0.0326398 & 0.0077994 & -0.434404 & 0.111103 \\
 \, P_2^{-1} & 27 & 0.0487394 & 1023 & 2.81247 & -0.33028 & 0.0036736 & -0.434123 & 0.207571 \\
\hline
\end{array}
\end{array}
$
\label{table:table294}
\end{table}
%
%
\begin{table}[ht]
\caption{Test Case $4$ ramification profile, $\mathscr{K}=132$}
$
\small
\begin{array}{c}
\begin{array}{|c|c|}
\hline
 \text{Singular point} & \text{Cycles} \\
\hline
\{\bar{s}\}_{127} &      (2,[33,1])\\
s_{\infty} & (5,2,[28,1])\\
\hline
\end{array}

		\end{array}
	$
\label{table:table55a}
	
\end{table}

%
%
%
%
\subsection{Test Case 5:  $50$-degree function expanded at the origin, $\mathscr{G}=2268$}
\label{testCase5}
\begin{equation}
  \begin{aligned}
    f_5(z,w)&=\left(2 z^6+\frac{1}{2} z^7-\frac{5}{4} z^{11}+4 z^{22}+\frac{29}{10} z^{34}-1 z^{40}-\frac{13}{2}
   z^{43}\right)+\left(\frac{3}{5} z^{10}+\frac{7}{4} z^{24}-\frac{1}{4} z^{50}\right) w^2\\
  &+\left(2 z^{17}+\frac{7}{2} z^{34}\right) w^3+\left(-\frac{3}{2} z^{30}+\frac{4}{3} z^{38}+\frac{8}{5} z^{42}\right) w^4+\left(-\frac{6}{5} z^2-\frac{1}{2} z^6+\frac{7}{3} z^{31}\right) w^9\\
  &+\left(-\frac{2}{5} z^{11}-\frac{3}{2} z^{26}+1 z^{45}\right)w^{10}+\left(\frac{7}{5} z^{24}-6 z^{32}-6 z^{49}\right) w^{14}\\
  &+\left(-\frac{3}{4} z^5+\frac{7}{3} z^{21}-\frac{1}{4} z^{26}+\frac{4}{5} z^{27}+\frac{4}{3} z^{32}-2 z^{36}+\frac{1}{3} z^{39}-\frac{3}{4} z^{41}-z^{43}\right)w^{16}\\
  &+\left(-6 z^{14}-2 z^{31}-z^{33}\right) w^{18}+\left(-2 z^{27}-\frac{8}{3} z^{50}\right) w^{22}+\left(4z^8+\frac{4}{5} z^{25}-\frac{3}{2} z^{27}\right) w^{24}\\
  &+\left(-3 z^4+\frac{8}{3} z^{22}-\frac{8}{5} z^{43}\right) w^{33}+\left(\frac{7}{3} z^{14}-\frac{3}{2} z^{18}\right) w^{34}+\left(-4+8 z^{13}-\frac{7}{4} z^{47}\right)w^{36}\\
  &+\left(z^2-\frac{1}{4} z^7\right) w^{38}+\left(-\frac{1}{2} z^{20}-z^{29}+z^{46}\right) w^{40}+\left(\frac{1}{3}z^{10}+\frac{7}{4} z^{11}+\frac{8}{5} z^{21}\right) w^{47}\\
  &+\left(\frac{2}{3} z^2+6 z^{26}+\frac{3}{5} z^{43}\right)w^{48}+\left(-z^9+\frac{1}{4} z^{13}+2 z^{14}+2 z^{18}+z^{36}-2 z^{44}\right) w^{49}\\
    &+\left(-\frac{1}{3}z^{23}-\frac{7}{2} z^{40}+z^{42}\right) w^{50}\\
  &=0
  \end{aligned}
	\label{eqn:equation2000}
	\end{equation}

 This function has $4584$ singular points.  The singular point at the origin was selected as $s_b$ and the Root Test indicated a comparison sequence of $s_2$ through $s_{100}$.  Table \ref{table:table477c} is the timing data, and Table \ref{table:table29a}, the Summary Report.
\begin{table}[h]
\caption{Test Case $5$ Timing Data}
$
\begin{array}{|c|c|c|c|c|c|}
\hline
 \begin{array}{c}
 \text{Singular} \\
  \text{points}
	\end{array}&
	\begin{array}{c}
 \text{Initial} \\
  \text{segments}
	\end{array}&
	\begin{array}{c}
 \text{Base} \\
  \text{expansion}
	\end{array}&
	\begin{array}{c}
 \text{Comparison} \\
  \text{expansions}
	\end{array}&
	\text{CT} &
	\text{IT} 
	\\
\hline
 $(4584,3.76\,h)$ & $(4584,49\,m)$ &$(6,50\,m)$ & $(100,11\,h)$ & $1.3\,m$ & $2.7\,m$  \\
 \hline 
\end{array} 
$
\label{table:table477c}
\end{table}
\begin{table} [ht]
\caption{Test Case $5$ Summary Report at the origin}
$
\small
\begin{array}{c}
\begin{array}{|c|c|c|c|c|c|c|c|c|}
\hline
 \text{Type} & \text{CLSP} & \text{R} & \text{Terms} & \text{a} & \text{b} & \text{c} & \text{d} & \text{Var} \\
\hline
 \, V_9^4 & 2 & 0.564982 & 1969 & 3.32621 & -0.245552 & 0.00915745 & -0.433441 & 0.165712 \\
 \, V_{27}^2 & 2 & 0.564982 & 1837 & 2.88746 & -0.0586786 & 0.0251435 & -0.435114 & 0.163804 \\
 \text{1$\, $}P_6^{-1} & 98 & 0.791246 & 2004 & 2.60055 & -0.314091 & 0.00555285 & -0.433077 & 0.200407 \\
 \text{2$\, $}P_6^{-1} & 94 & 0.786365 & 2004 & 2.47551 & -0.340257 & 0.0073696 & -0.432927 & 0.215718 \\
 L_{-7} & 92 & 0.75919 & 1020 & 0.865693 & -0.405906 & 0.00141605 & -0.43457 & 0.285315 \\
 L_{-14} & 92 & 0.75919 & 1017 & 0.694339 & -0.323955 & 0.00148271 & -0.434907 & 0.191351 \\
\hline
\end{array}
\end{array}
$
\label{table:table29a}
\end{table}
%
%

\begin{table}[ht]
\caption{Test Case $5$ Ramification profile, $\mathscr{K}=4634$}
$
\small
\begin{array}{c}
\begin{array}{|c|c|}
\hline
 \text{Singular point} & \text{Cycles} \\
\hline
 s_1 &(9,27,[2,6],[2,1])\\
\{p_i\}_{19} & [50,1] \\
\{\bar{s}\}_{4564} &(2,[48,1])\\
s_{\infty} & (14,13,2,[21,1]))\\
\hline
\end{array}
\\
\scriptsize \textbf{$\{p_i\}_{19}=\{1,229,230,233,234,238,239,248,249,262,263,276,277,300,301,324,325,334,4537,4538\}$} \\
		\end{array}
	$
\label{table:table55b}
\end{table}
%
%
%
%
%
\subsection{Test Case $6$: $25$ degree function with complex coefficients expanded at $s_{29}$, $\mathscr{G}=326$}

\label{testCase6}

\begin{equation}
\begin{aligned}
f_6(z,w)&=-\left(\frac{311}{20}i+\frac{467}{30}\right)-\left(\frac{16}{3}i+\frac{3}{2}\right) z+\left(\frac{1}{4}-\frac{3}{4}i\right)z^3+\left(\frac{9}{4}i+3\right) z^4-\left(\frac{3}{10}-\frac{21}{5}i\right) z^5\\
	&-\left(1-\frac{7}{3}i\right)
   z^7-\frac{5}{3}z^8-\left(\frac{12}{5}i+1\right) z^9-\left(2 i+\frac{5}{4}\right) z^{11}-\left(\frac{1}{4}-3 i\right)z^{12}\\
	&-\left(\frac{27}{5}i+3\right) z^{13}+\left(\frac{19}{6}i+7\right) z^{14}+(3 i-2)   z^{15}\\
	&+\biggr(-\left(\frac{59}{60}i+\frac{13}{10}\right)+(1-3 i) z^{12}+\left(4 i+\frac{4}{5}\right)z^{13}\biggr)w\\
	&+\biggr(-\left(\frac{71}{10}-\frac{1}{6}i\right)-\left(\frac{8}{3}-6 i\right) z^2+(-i-8)
   z^8\biggr)w^6\\
	&+\biggr(\frac{15}{2}i+\frac{116}{15})w^8+(\left(\frac{12}{5}-\frac{15}{4}i\right)+\frac{2 i z^3}{5}+\left(-7
   i-\frac{5}{2}\right)z^{10}\biggr)w^9\\
	&+\biggr(\frac{13}{2}-\frac{32}{5}i\biggr)w^{10}\\
	 &+\biggr(-\left(\frac{19}{4}i+\frac{5}{3}\right)-\left(3-\frac{6}{5}i\right)
   z^3+\left(\frac{3}{2}i-\frac{3}{4}\right)z^{13}\biggr)w^{13}\\
	&+\biggr(\left(\frac{47}{15}i+\frac{155}{12}\right)-\left(\frac{1}{6}-\frac{1}{2}i\right)
   z^2-\left(\frac{1}{5}-\frac{5}{3}i\right) z^{10}+\left(\frac{3}{2}i-2\right)
   z^{14}\biggr)w^{14}\\
	&+\biggr(-\left(\frac{109}{12}-\frac{37}{12}i\right)-\frac{i z^5}{2}\biggr)w^{18}\\
	&+\biggr(\left(5i+\frac{6}{5}\right)-\left(\frac{1}{5}i+\frac{7}{4}\right) z^{12}+\left(\frac{8}{5}i+7\right)
   z^{14}\biggr)w^{21}\\
	& +\biggr(\left(\frac{2}{5}-\frac{13}{10}i\right)-\frac{8}{3}z^7+(-2 i-1) z^{10}\biggr)w^{25}
		\end{aligned}
		\end{equation}

 This function has $660$ finite singular points and was expanded at pole $s_{29}$.  Table \ref{table:table47c} gives the timing data and Tables \ref{table:table29} and \ref{table:table56b} detail the accuracy and ramification data.

%
%
\begin{table}[h]
\caption{Test Case $6$ Timing Data at $s_{29}$}
$
\begin{array}{|c|c|c|c|c|c|}
\hline
 \begin{array}{c}
 \text{Singular} \\
  \text{points}
	\end{array}&
	\begin{array}{c}
 \text{Initial} \\
  \text{segments}
	\end{array}&
	\begin{array}{c}
 \text{Base gen.} \\
  \text{expansions}
	\end{array}&
	\begin{array}{c}
 \text{Comparison} \\
  \text{expansions}
	\end{array}&
	\text{CT} &
	\text{IT} 
	\\
\hline
 $(660,30\,m)$ & $(660,3.3\,m)$ &$(22,2\,h)$ & $(27,29\,m)$ & $1.3\,m$ & $3.4\,m$  \\
 \hline 
\end{array} 
$
\label{table:table47c}
\end{table}		
\begin{table} [ht]
\caption{Test Case $6$ Summary Report at $s_{29}$}
$
\small
\begin{array}{l}
\begin{array}{|c|c|c|c|c|c|c|c|c|}
\hline
 \text{Type} & \text{CLSP} & \text{R} & \text{Terms} & \text{a} & \text{b} & \text{c} & \text{d} & \text{Var} \\
\hline
 \text{1T} & 38 & 0.0167133 & 1024 & 3.17032 & -0.41029 & 0.00192871 & -0.434412 & 0.219225 \\
 \text{2T} & 56 & 0.0982504 & 1024 & 3.69649 & -0.368397 & 0.00189034 & -0.434359 & 0.242937 \\
 \text{3T} & 177 & 0.173224 & 1024 & 3.68729 & -0.446054 & 0.00194646 & -0.434356 & 0.248846 \\
 \text{4T} & 354 & 0.208561 & 1024 & 3.54559 & -0.371434 & 0.00189941 & -0.434436 & 0.254165 \\
 \text{5T} & 10 & 0.197401 & 1024 & 3.98117 & -0.433098 & 0.0018211 & -0.434431 & 0.215548 \\
 \text{6T} & 232 & 0.156337 & 1024 & 3.64785 & -0.336867 & 0.00194168 & -0.434384 & 0.227871 \\
 \text{7T} & 108 & 0.149061 & 1024 & 3.68843 & -0.385887 & 0.00197982 & -0.434356 & 0.249595 \\
 \text{8T} & 41 & 0.0672525 & 1024 & 3.48234 & -0.389387 & 0.00182101 & -0.434454 & 0.241994 \\
 \text{9T} & 125 & 0.101533 & 1024 & 3.82856 & -0.431569 & 0.00196967 & -0.434425 & 0.258889 \\
 \text{10T} & 13 & 0.0970926 & 1024 & 3.53643 & -0.479354 & 0.00191281 & -0.43444 & 0.25699 \\
 \text{11T} & 125 & 0.101533 & 1024 & 3.90238 & -0.37848 & 0.00192665 & -0.434415 & 0.25622 \\
 \text{12T} & 264 & 0.16936 & 1024 & 3.74861 & -0.428541 & 0.00204153 & -0.434296 & 0.243689 \\
 \text{13T} & 20 & 0.0362971 & 1024 & 3.44335 & -0.374599 & 0.0018982 & -0.434408 & 0.245484 \\
 \text{14T} & 129 & 0.109811 & 1024 & 3.58688 & -0.357851 & 0.00184152 & -0.434469 & 0.225132 \\
 \text{15T} & 195 & 0.140624 & 1024 & 3.62277 & -0.324291 & 0.00190887 & -0.434374 & 0.228172 \\
 \text{16T} & 86 & 0.138713 & 1024 & 3.56908 & -0.453317 & 0.00192339 & -0.434365 & 0.252658 \\
 \text{17T} & 210 & 0.142419 & 1024 & 3.74149 & -0.386165 & 0.00196671 & -0.43432 & 0.233369 \\
 \text{18T} & 181 & 0.201697 & 1024 & 3.64415 & -0.423745 & 0.00191286 & -0.43442 & 0.276286 \\
 \text{19T} & 120 & 0.16032 & 1024 & 3.32947 & -0.453331 & 0.0020065 & -0.434385 & 0.269883 \\
 \text{20T} & 173 & 0.175719 & 1024 & 3.62089 & -0.440735 & 0.00191062 & -0.434385 & 0.24661 \\
 \text{21T} & 35 & 0.0170704 & 1024 & 3.26908 & -0.40917 & 0.00189352 & -0.434389 & 0.238545 \\
 \, P_4^{-1} & 38 & 0.0167133 & 1021 & 3.08615 & -0.0848102 & 0.00655571 & -0.434368 & 0.160934 \\
\hline
\end{array}
\\
 \scriptsize \textbf{In this case $R=|s_{29}-s_{CLSP}|$}\\
 
		\end{array}
$
\label{table:table29}
\end{table}
%
%

\begin{table}[ht]
\caption{Test Case $6$ Ramification profile, $\mathscr{K}=700$}
$
\small
\begin{array}{c}
\begin{array}{|c|c|}
\hline
 \text{Singular point} & \text{Cycles} \\
\hline
 \{p_i\}_{10} &(4,[21,1])\\
\{\bar{s}\}_{650} &      (2,[23,1])\\
s_{\infty} & (21,[4,1])\\
\hline
\end{array}
\\
 \scriptsize \textbf{$\{p_i\}=\{29,30,37,43,76,88,212,316,486,509\}$}\\
 
		\end{array}
	$
\label{table:table56b}
\end{table}
\section{Conclusions}
\begin{enumerate}
\item A numerical approach to Newton polygon initially seems problematic.  However this work includes several error checking algorithms:
\begin{enumerate}
	\item \textbf{Cycle check-sum: } The sum of the conjugate classes at any expansion center must equal to the degree of the function in $w$,
	\item \textbf{Global cycle check-sum: } The Riemann-Hurwitz sum must be a positive number,
	\item \textbf{Removal of coefficient zeros: } In order to minimize residual errors, coefficient zeros are removed from $f$ and its iterates before processing by the Newton polygon algorithm,
	\item \textbf{Precision monitoring: } The precision of the calculations are monitored throughout the algorithms and terminate the analysis if it drops below $900$ digits,
	\item \textbf{Accuracy results: } The accuracy results of a set of expansions would not follow the log-linear trend with low variance if a branch was incorrectly computed.
\end{enumerate}
 
These measures reduce the potential of errors.   However, there exists functions which can usurp this numeric approach.  These would include the following scenarios:

\begin{enumerate}
  	\item \textbf{Singular size limits: } Functions with very small singular sizes coupled with very large exponents sufficient to compromise a realistic level of precision achievable in a reasonable amount of time,
	\item \textbf{High polygon iterates: } Functions which entail multiple polygon iterations sufficient to reduce the precision of the characteristic equations below a reasonable level of numeric precision,
	\item \textbf{High polygon iterates: } It is likely there are functions with arbitrary polygon iterates which would decrease the precision of the calculations beyond any effort to keep the results above a minimum level,  
	\item \textbf{High poly mod function: } Functions with extremely large gaps between successive expansion terms would cause the modular operation of (\ref{eqn:eqn989}) to exceed the (arbitrary) maximum number of zero term iterations.  In this case, an infinite power expansion would be identified as a fractional polynomial.  
	\end{enumerate}

  \item Identifying conjugate classes and only expanding generator series of each class is an improvement to the standard approach to expanding all initial segments of a Newton Polygon expansion.  The greatest time saving is when the conjugate set at an expansion center is highly ramified.  An example of this is the $50$-degree function of test case $4$ which only required iterating six generator series in $50$ minutes rather than an estimated seven hours to generate the full $50$ set. 
  \item  The test cases were designed to stress-test the series comparison and integration tests with complex functions.  With minor tuning of the methods, CLSPs and corresponding radii of convergence results agreed well with the estimated values determined by the Root Test when the separation tolerance was set to $1/10$ the separation distance of the branch values.  In cases where the comparison test and integration test succeeded in identifying CLSPs, both identified the same CLSP for each branch.  In Test Case $3$ where the comparison test failed, this was due to extremely small singular perimeters on the order of $10^{-35}$ causing the base expansions to be evaluated extremely close to their radii of convergence.  However the integration test succeeded with proper adjustments to the integration method.
		\item The accuracy and order functions were found to agree well with the actual accuracies of the series as shown by the accompanying low variances of the fit function showing $A(r_f,o)$ to be a robust predictor of accuracy in the testing range.
	\item The close approximation of the Root Test to the CLSPs shows the Root Test to be a reliable means of estimating radii of convergences.
	\item  This work opens the subject to further research:
	\begin{enumerate}
		\item Find and analyze function which entail more polygon iterations,
		\item Find and analyze functions with fractional polynomial solutions,
		\item Further fine-tuning the algorithm as problem cases arise,
		\item Reducing execution time,
		\item Identifying problem functions and updating the method to accommodate them. 
	\end{enumerate}
	
	\end{enumerate}

\section{}
\label{appendix:appA}
\begin{center}
Appendix A: Branch types used in this paper
\end{center}

In the following branch descriptors, all exponents $\displaystyle \frac{q}{p}$ of a series are presumed placed under a least common denominator $p$.

\begin{description}
\item[\textbf{Type \boldsymbol{$T$}}]  Power series with positive integer powers (Taylor series).  These are $1$-cycle branches.
\vspace{5pt}
 \item[\textbf{Type \boldsymbol{$E$}}] $1$-cycle $T$ branch with a removable singular point at its center.
\vspace{5pt}
\item[\textbf{Type \boldsymbol{$F_p^q$}}] $p$-cycle branch with $p>1$ of order $q$ with non-negative exponents and lowest non-zero exponent $\displaystyle \frac{q}{p}$ with $q>p$.  These branches are multi-valued consisting of $p$ single-valued sheets with a finite tangent at the singular point.  An example $F_2^3$ series is $z^{3/2}+z^2+\cdots$.
\vspace{5pt}
\item[\textbf{Type \boldsymbol{$V_p^q$}}] $p$-cycle branch with $p>1$ of order $q$ with non-negative exponents and lowest non-zero exponent $\displaystyle \frac{q}{p}$ with $q<p$ and vertical tangent at center of expansion.  An example $V_4^3$ series is $z^{3/4}+z^2+\cdots$.
\vspace{5pt}
\item[\textbf{Type \boldsymbol{$P_p^q$}}]  $p$-cycle branch unbounded at center with $p>1$ of order $q$ having negative exponents with lowest negative exponent $\displaystyle\frac{q}{p}$.   An example $3$-cycle $P$ series of order $-1$ is $z^{-1/3}+z^2+\cdots$.  An example $3$-cycle $P$ series of order $-5$ is $z^{-5/3}+z^{-1/3}+\cdots$. 
 \vspace{5pt}
\item[\textbf{Type \boldsymbol{$L^q$}}]  Branch with Laurent series of order $q$ as the Puiseux series.  An example $L^{-2}$ series is $1/z^2+1/z+z^2+\cdots$
\end{description}


\begin{thebibliography}{10}  

\bibitem {Bliss} Bliss, Gilbert A. \textit{Algebraic Functions}. New York: Dover Publications, Inc., 2004.

\bibitem {Brown} Brown, James and Ruel Churchill. \textit{Complex Variables and Applications}. New York:  McGraw Hill, 2004

\bibitem {Chud} Chudnovsky, D.V. and G.V. Chudnovsky. ``On Expansion of Algebraic Functions in Power and Puiseux Series''. Journal of Complexity \textbf{2}, 271-294 (1986).

\bibitem {Kung} Kung, H.T. and J. Traub, ``All Algebraic Functions can be Computed Fast''. J. Assoc. Comput. Mach. \textbf{25}, 245-260.

\bibitem {Marku} Markushevich,A.I.,1967.\textit{Theory of Functions of a Complex Variable.Vol.III.} PrenticeHall,
Englewood Cli?s, N. J.

\bibitem {Marsden} Marsden, Jerrold and Michael Hoffman. \textit{Basic Complex Analysis}. New York: W.H Freeman and Company, 1999.
\bibitem {Milioto} Milioto,Dominic C. (2018, Dec. 8). Algebraic Functions, Retrieved from \href{https://jujusdiaries.com}{Algebraic Functions and Iterated Exponentials}.
\bibitem {Milioto2} Milioto,Dominic C. (2018,Jan.13).  On the branching geometry of algebraic functions, Retrieved from \href{https://arxiv.org/abs/1901.03996}{On the branching geometry of algebraic functions}
\bibitem {Milioto3} Milioto,Dominic C. (2021,Dec.13). Determining radii of convergence of fractional power expansions around singular points of algebraic functions, Retrieved from \href{https://arxiv.org/pdf/2111.11883.pdf}{Determining radii of convergence of fractional power expansions around singular points of algebraic functions}

\bibitem {Nowak} Nowak, Krzysztof.\ \textit{Some Elementary Proofs of Puiseux's Theorem}.  Universitatis Iagellonicae ACTA Mathematica, Fasciculus XXXVIII, 2000. 
\bibitem {Walker} Walker, Robert J.  \textit{Algebraic Curves}. Princeton:  Princeton University Press, 1956.

\bibitem {Willis} Willis, Nicholas J., Didier, Annie K., Sonnanburg, Kevin M. \textit{How to Compute a Puiseux Expansion},
arXiv: 0807.4674.1 [math.AG] 29 July, 2008
\end{thebibliography}
\end{document}